\documentclass[]{AntiMath24}
\usepackage[OT2,T1]{fontenc}
\usepackage[utf8]{inputenc}
\usepackage{enumerate}
\usepackage{soul}

\input cyracc.def
\newcommand\cyrfamily{\fontencoding{OT2}\fontfamily{wncyr}
\selectfont\cyracc}
\DeclareTextFontCommand{\textcyr}{\cyrfamily}

\usepackage[english,russian,german,polish]{babel} 
\usepackage{lastpage} 
\usepackage{comment}

\usepackage{etex} 
\usepackage{filecontents}

\usepackage{multirow}

\usepackage{amsmath}
\usepackage{color}

\usepackage[colorlinks=true,linkcolor=blue]{hyperref}
  \hypersetup{
    pdftitle={Barbara Rokowska's (1926--2012) combinatorial research}, 
    pdfauthor={Krzysztof J. Szajowski}, 
    bookmarks=true,
    colorlinks,
    urlcolor=blue,
    filecolor=magenta,
    citecolor=green, 
    linkbordercolor={1 1 1}, 
    citebordercolor={1 1 1},  
    urlbordercolor={ 1 1 1}  
  } 
\RequirePackage[numbers]{natbib}

\RequirePackage[hyperpageref=true]{backref}

    \renewcommand*{\backref}[1]{}  
    \renewcommand*{\backrefalt}[4]{
       \ifcase #1 
          Nie cytowana.
       \or
          Cyt. na str. #2.
       \else
          Cyt. na str. #2.
       \fi}

\def\BD{\textbf{BD}}
\def\BIBD{\textbf{BIBD}}
\def\CM{\textit{CM}}
\def\SD{\textbf{SD}}
\def\SS{\textbf{SS}}
\def\STS{\textbf{STS}}
\def\UiPWr{\textbf{UiPWr}}
\def\UWr{\textbf{UWr}}
\def\PWr{\textbf{PWr}}

\def\IM{\textbf{IM}}

\def\PAN{\textbf{PAN}}

\def\IMiFT{\textbf{IMiFT}}
\def\cB{\mathcal{B}}
\newcommand*{\MRid}[1]{\href{http://www.ams.org/mathscinet/search/author.html?mrauthid=#1}{MRid:#1}} 

\pdfoutput=1 
\renewcommand{\contentsname}{}
\usepackage{graphicx}
\usepackage{wrapfig}
\usepackage{subcaption}
\graphicspath{{Pictures/},{./Figures/},{../Figures/},{./Pictures/}}
\renewcommand{\figurename}{Fig.}
\usepackage{float}
\firstpage{1}
\LogoG{\includegraphics[width=0.18\textwidth]{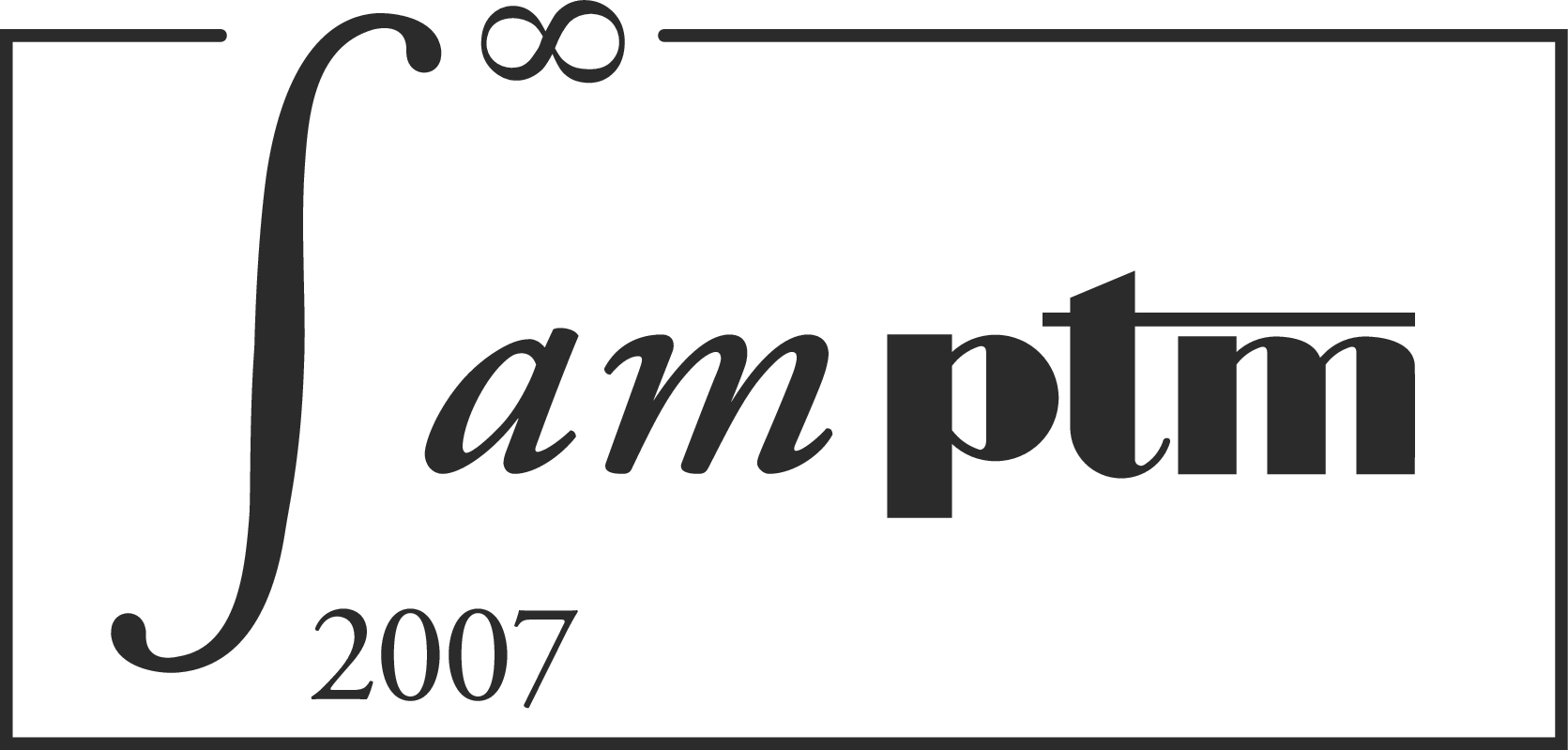}}
\volume{18}
\fasc{1}
\years{2024}
\LogoGMS{\includegraphics[width=0.18\textwidth]{am.png}}
\volumeMS{18}
\numberMS{18}
\wwwtrue
\secnameMS{HISTORY OF APPLIED MATHEMATICS}
\pages{1–31}
\received[31st of December 2024]{30th of September 2023}
\lastrevision{23rd August 2024}
\logo{}{}{}
\def\doinum{10.14708/am.v18i0.71xx}

\title[Barbara Rokowska's (1926--2012)] {Barbara Rokowska's combinatorial research\\with her extensive biography (1926--2012)}

\tpauthortrue 

\author[K.J. Szajowski]{Krzysztof J. Szajowski\orcid{0000-0001-9834-9929}}
\affiliation{\vspace{2ex}Wroclaw University of Science and Technology}
\affiliation{Faculty of Pure and Applied Mathematics}
\address{Wybrzeże Wyspiańskiego 27, 50-370 Wroclaw, Poland}
\email{Krzysztof.Szajowski@pwr.edu.pl}
\city{Wrocław}


\keywordsPL{Systemy Steinera; system bloków; liczby pierwsze; ciała Galois; konstrukcje kombinatoryczne}

\comm{Stanisław Domoradzki} 
\begin{document}
\vspace{-6ex}
\setcounter{page}{1} 
\selectlanguage{english}
\Polskifalse\renewcommand{\figurename}{Fig.}
\begin{abstract}
We discuss the significance of some interesting results by Barbara Rokowska about combinatorial constructions. Her interest in finite mathematics and number theory began with an embellishment and detailing of some work by Erdős. Rokowska and Schinzel then solved the problem posed by Paul Erdős concerning the existence of prime numbers of a certain kind. Her subsequent work highlighted the difficulty in constructing Steiner systems ({\SS}s)\label{KSzSS} with certain properties and showed the importance of rigorous proof techniques in this area of mathematics. This is the first such summary of the main results obtained by Rokowska, her collaborators and PhD students.

A biography of Barbara Rokowska has been added as an appendix.
\end{abstract}
\vspace{-5ex}
\abbreviations{\label{skroty}W artykule oznaczono skrótem głównie instytucje, czasem terminy z zakresu badań.\\
 
\begin{tabular}{@{}llll}
\BIBD & a balanced incomplete block design (v.~\pageref{BRbibd}, App.~\ref{BRappBIBD})\\
\SS    & Steiner triple system (v.~\pageref{KSzSS})\\
\STS    & Steiner triple system (v.~\pageref{KSzSTS})\\
\UiPWr   & Uniwersytet i Politechnika Wrocławska\\
&\textsl{ie University and Technical University of Wrocław} (v. \pageref{BRUiPWr}) \\
\UWr & Uniwersytet Wrocławski (\textsl{ie Wrocław University} v. \pageref{BR_UWr}) \\
\PAN & Polish Accademy of Science (v. \pageref{BR_PAN}) \\
\PWr  & Politechnika Wrocławska (\textsl{ie Wrocław Technical University } v.~\pageref{BRPWr})\\
\IM & Mathematical Institute (v. \pageref{BR_IM}) \\
\IMiFT & Instytut Matematyki i Fizyki Teoretycznej\\
&\textsl{ie Institute of Mathematics and Theoretical Physics}  (v. \pageref{BR_IMiFT})
\end{tabular}}
\medskip 

\section{Introduction.} 
The purpose of this note is to discuss the achievements in mathematics of Barbara Rokowska, a professor at Wroclaw University of Technology, whose abbreviated biography is included in the supplementary materials. Rokowska is one of the few examples known to me that one's path in life should be chosen prudently, and when seeing a wrong choice, one should correct it--even when it means incurring significant costs and sacrifices. With some delay in her primary education caused by World War II, and passing her high school diploma, she chose to study at the Faculty of Polish Philology at the University of Wrocław\footnote{At that time it was administered as University and Technical High School of Wroclaw (\UiPWr\label{BRUiPWr} v.~\citeauthor{Chmielewski2007:WSA}~ (\citeyear{Chmielewski2007:WSA})).}, which she graduated on time in 1951 at the age of 25. After graduation, she took up a job, but conversations with friends led her to believe that the education she had received and the profession associated with it were not the ones that would allow her to pursue her interests and passions. 

\begin{figure}[H]
    \centering
    \begin{subfigure}[b]{0.494\textwidth}
    \includegraphics[width=\textwidth]{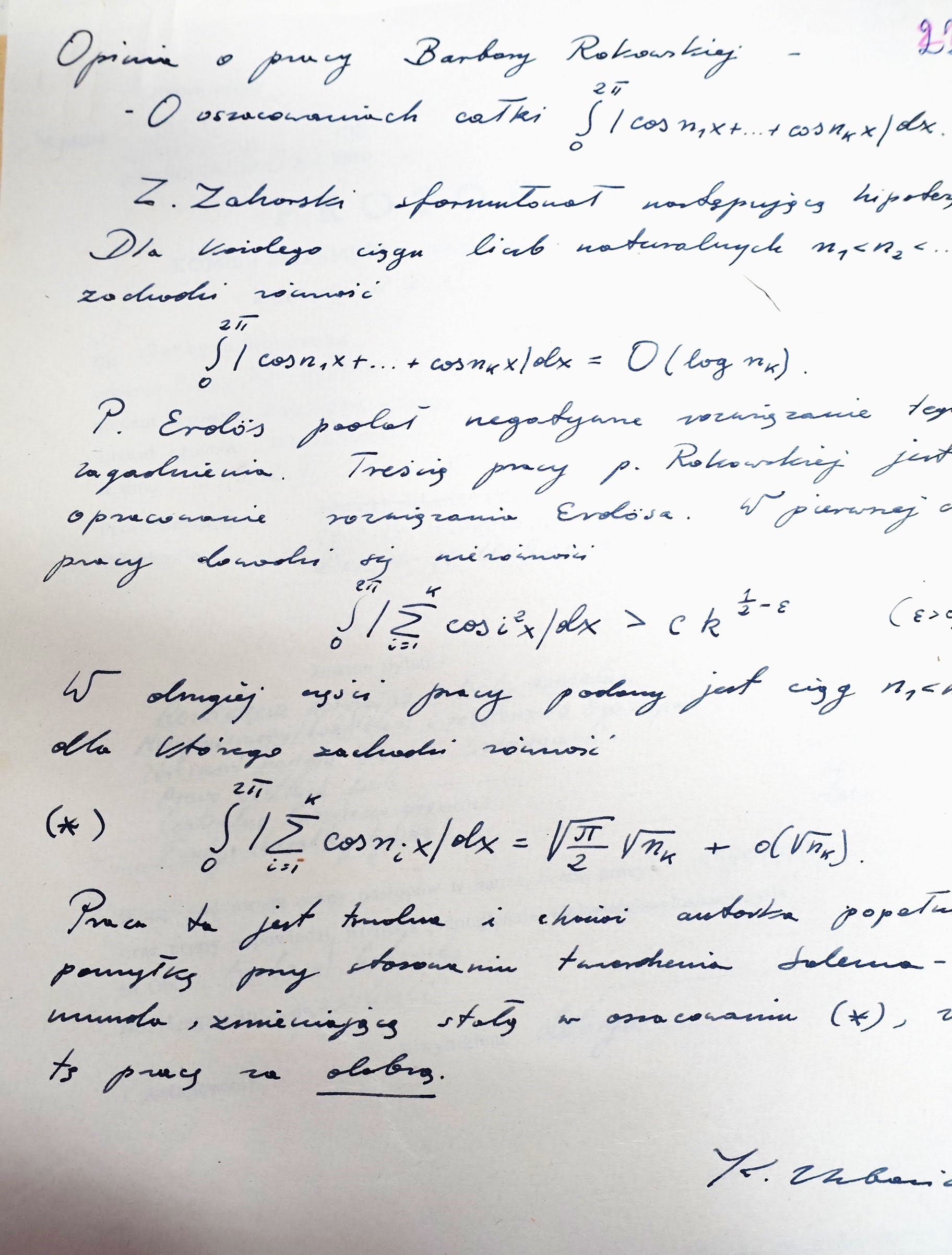}
       \caption{\small\label{BRRevKUMgr} Urbanik's review of the thesis.}
    \end{subfigure}
    \hfill
    \begin{subfigure}[b]{0.494\textwidth}
        \includegraphics[width=0.962\textwidth,height=0.4\textheight]{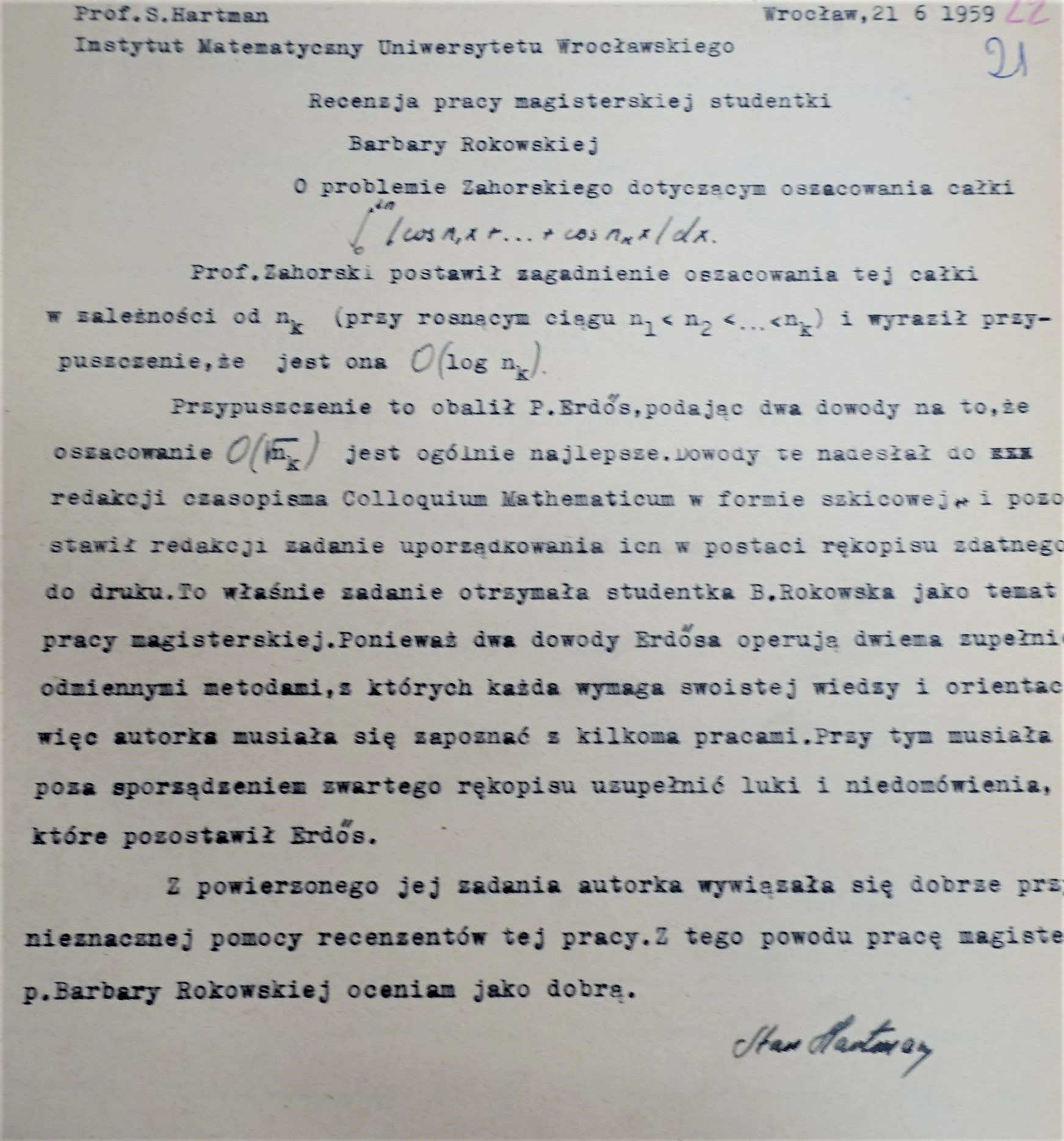}
      \caption{\small\label{BRRevSHMgr} Hartman's review of the thesis.}
    \end{subfigure}
    \caption{\label{BR2RevMgr} Master thesis records (v. \cite[mks 23,24 AUWr]{BRokowska1954:UWr})}
~ 
\end{figure}
In 1954, she takes the exam to study full-time mathematics at the Faculty of Mathematics, Physics and Chemistry at the University of Wrocław (\UWr\label{BR_UWr}). Initially, she combines her education in the humanities with her studies by conducting research in the history of mathematics. Thus, originally, the reason for studying mathematics, as she herself admits, was her interest in the history of mathematics (v. p.~32--33 in~\citeauthor{Rok1996:Jezus}~(\citeyear{Rok1996:Jezus})). 

In the initial period of mathematical studies, she supported herself by tutoring (private teaching).  As she achieved results, she was awarded a scholarship that allowed her to spend less time on tutoring. At the same time, she started working as a technical editor for the publications of the State Mathematical Institute in Wrocław. One of the journals published by the PIM in Wrocław was Colloquim Mathematicum (further~\CM), the editors of which were lecturers in the mathematics department of the \UWr, where she was studying. At the time when students were receiving thesis topics, the \CM\ received a paper from the thesis \citeauthor{Erdos1960:Zahorski}~(\citeyear{Erdos1960:Zahorski}) which, although very interesting, was incomplete and poorly edited. The task of correcting the inaccuracies was given to Rokowska. The task was completed with a good grade (see Fig.~\ref{BRRevKUMgr} and Fig.~\ref{BRRevSHMgr}) and Barbara Rokowska received her Master's degree in Mathematics~(see ~\cite[mks 23,24 AUWr]{BRokowska1954:UWr}, Fig.~\ref{BRMgr}). More details are given in the addendum, chapter~\ref{KSzBRBio}.

Thus, a natural area of interest after graduation, given the topic of her master's thesis, is number theory. The first publication concerns the periodicity of certain number sequences, a problem posed by Sierpinski. The result was published by Rokowska in \textsl{Wiadomosci Matematyczne} (v.~ [BR\ref{Rok1959:First}]). Together with Schinzel in [BR\ref{RokSchi1960:Erdos}], she analyzed the problem posed by~\citeauthor{Erd1960:MR}~(\citeyear{Erd1960:MR})\footnote{Does there exist a prime $p>5$ for which all the numbers $2!, 3!, \ldots, (p-1)!$ are incongruent $\mod p$?} giving a partial answer. She then turns her attention to problems in which the importance of number theory is significant - finite geometries and systems of blocks (v.~mps sheet 29 in \citeauthor{BRok1966:PhD}~(\citeyear{BRok1966:PhD})). The main results of the doctoral dissertation were published on the paper [BR\ref{Rok1967:PhD1}]. More about the theses of the doctorate is in~chapter~\ref{sBRPhD}. 
\begin{figure}[th!]
\centerline{\includegraphics[width=1.1\textwidth]{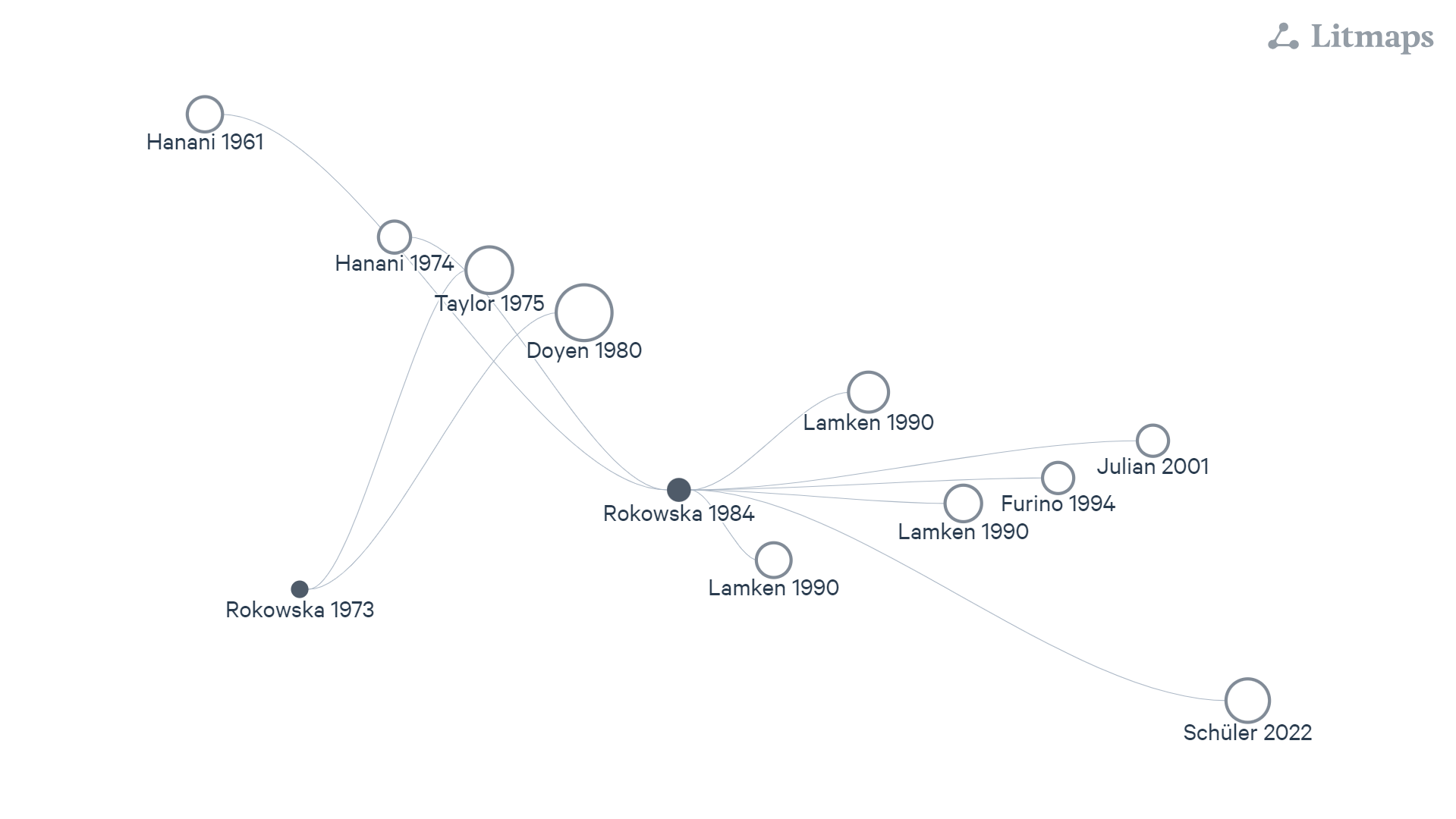}}
\caption{\label{BR_LipMap1}Links to the papers~[BR\ref{Rok1973:STS73}] and [BR\ref{Rok1984:STS84}] from related research articles.}
\end{figure}
At that time, she participated in the seminar taught by Associate Professor Lucjan Szamkolowicz~(1927--1984). The close relations in the scientific discussions about combinatorial configurations at that time were held with Mieczysław Wojtas and Krystyna Pukanow. It should be emphasized that soon after receiving her Ph.D., there were organizational changes at her workplace. In 1968, the Faculty of Fundamental Problems of Technology was established with the Institute of Physics and the Institute of Mathematics and Theoretical Physics. Almost all mathematicians employed at TU were associated with this Faculty. There was an institutional emphasis on staff research. One of the seminars was devoted to combinatorial analysis and was conducted by Associate Professor Zbigniew Romanowicz, Barbara Rokowska, and Lucjan Szamkowicz. Dr. Rokowska focused her attention at that time on the construction of block systems that generate non-isomorphic structures. The basic concepts of this area of combinatorics, definitions of the terms appearing in Rokowska's work, are contained in the supplementary chapter~\ref{BRappBIBD}. 

Before the end of 1974, she submitted her habilitation thesis to the Faculty of Mathematics, Physics, and Chemistry at the University of Wrocław. The Faculty Council implemented the procedure by appointing a three-member committee. The committee's task was to assess the completeness of the dossier and propose three reviewers. The dissertation theses collected results published or submitted for publication on possibly significantly different structures (v.~Section~\ref{BRHab} for the description of the results included in the habilitation dissertation).  

The procedure, awaiting reviews and confirmation of the degree by the State Commission for Academic Degrees and Titles, continued until the end of 1977. A discussion of the various stages of the proceedings is included in Section~\ref{KSzBRBio}. Dr. Rokowska's description of the habilitation procedure was typical of the period. It is included here to perpetuate the atmosphere around this procedure and also to show the great caution of mathematicians in expanding the field of mathematical research. 

A separate chapter discusses the subject of doctorates defended with Dr. Rokowska's supervision. The study of this issue should be expanded to include the relevance of the topics of these dissertations in relation to current research in combinatorial analysis. 

A set of additional chapters includes a biography of Dr. Rokowska and a complete list of her publications. To show the level of analysis of the dissertation, one of the reviews of the postdoctoral dissertation is included in the Appendix~\ref{BRhabRecSH}. 

\section{\label{sBRPhD} Dissertation on a Steiner systems.} She prepared her doctoral thesis {\color{blue}{Certain new constructions of Steiner systems}} under the supervision of Professor Czesław Ryll-Nardzewski (v. \citeauthor{BRok1966:PhD}~(\citeyear{BRok1966:PhD})). 
Her PhD defense took place on March 9, 1966, before the Scientific Council of the Department of Mathematics, Physics, and Chemistry\label{BR_UWr}\footnote{Her first paper was noted by Mathematical Review in 1959, giving her an identifier in that database \MRid{149940}. Related to this identifier is a record in the database \href{https://www.mathgenealogy.org/id.php?id=231649}{Mathematical Genealogy Project}, where her doctoral students are also listed, as well as the titles of their dissertations and the year of their defense}. A formal description of the procedure can be found in Section~\ref{BRPhDproc}. The results of the PhD dissertation can be traced in [BR\ref{Rok1967:PhD1}]. 

Here, to describe the achievements contained in the doctorate, we will use quotes from the review (cf. \citeauthor{Rokowska1966}~(\citeyear{Rokowska1966})). 
\begin{quote}[Doc. Dr Jan Mycielski]
The author's result is now the construction $S(v,4,1)=(V,\cB)$, where $|V|=v$ is appropriate, e.g. $v=mn\pm 1$, where $m,n\equiv 1 \vee 3\mod 6$. They represent a continuation of the research of Professor H.~Hanani of Haifa. 

Combinatorics in general and the topics of this work in particular have a special charm. It is characterized by the maximum simplicity of its basic concepts, since it basically operates only on finite sets. Also, its methods are usually modest, including some elements of number theory and the theory of finite groups and finite rings. Therefore, creativity in this field is very difficult. Its only basis is combinatorial fantasy, the kind of fantasy that produces the most surprising and valuable results in mathematics.

The above remarks testify well to the abilities of the Author, who has already had several other works in number theory, and one hopes that she will continue to work creatively. 
\end{quote}

\begin{quote}[Doc. Dr Andrzej Schinzel] The Steiner systems referred to in the paper are families of four-element subsets of a certain finite set $V$ with the property that every three-element subset of this $V$ is contained in one of the elements of this family. The author shows how, having suitable Steiner systems for the $m+1$ element set and the $n+1$ element set, to construct a Steiner system for the $nm+1$ element set, and gives three other related constructions.

The devising of these constructions undoubtedly required considerable ingenuity and was based on an analysis of Hannanie's existing work on the subject.
\end{quote}

The thesis was published in [BR\ref{Rok1967:PhD1}]. \citeauthor{Hanani1967:MR}~(\citeyear{Hanani1967:MR}) wrote about this result following:
\begin{quote} 
Let $E$ be a set of $v$ elements. By a quadruple system is understood a family of subsets of $E$ having four elements each, such that every triple of elements of $E$ is included in exactly one of these subsets. The reviewer proved that a necessary and sufficient condition for the existence of a quadruple system is that $v\equiv 2\text{ or }4(\hspace{-.5em}\mod 6\; )$ \cite{Hanani1960:Quadruple,Hanani1963:Tact}\footnote{It is worth noting at this point that the subject classification in the area of discrete mathematics changed in 1972/72. Dr. Rokowska's work had an assignment of \textsl{05.20 (1959-1972) Block designs}, while from 1973 it will rather be classified under \textsl{05B05 (1973-now) Combinatorial aspects of block designs}. The changes applied during this period indicate a significant number of works in this subject area and clearly distinguished additional subclasses.}.
The author gives some new constructions of quadruple systems not isomorphic to the known ones.
\end{quote}
The laconic nature of this opinion and the one-sentence conclusion, with no further commentary, is a token of appreciation for the result Rokowska discussed. 
\begin{figure}[th!]
    \centering
    \begin{subfigure}[b]{0.492\textwidth}
    \includegraphics[width=0.96\textwidth]{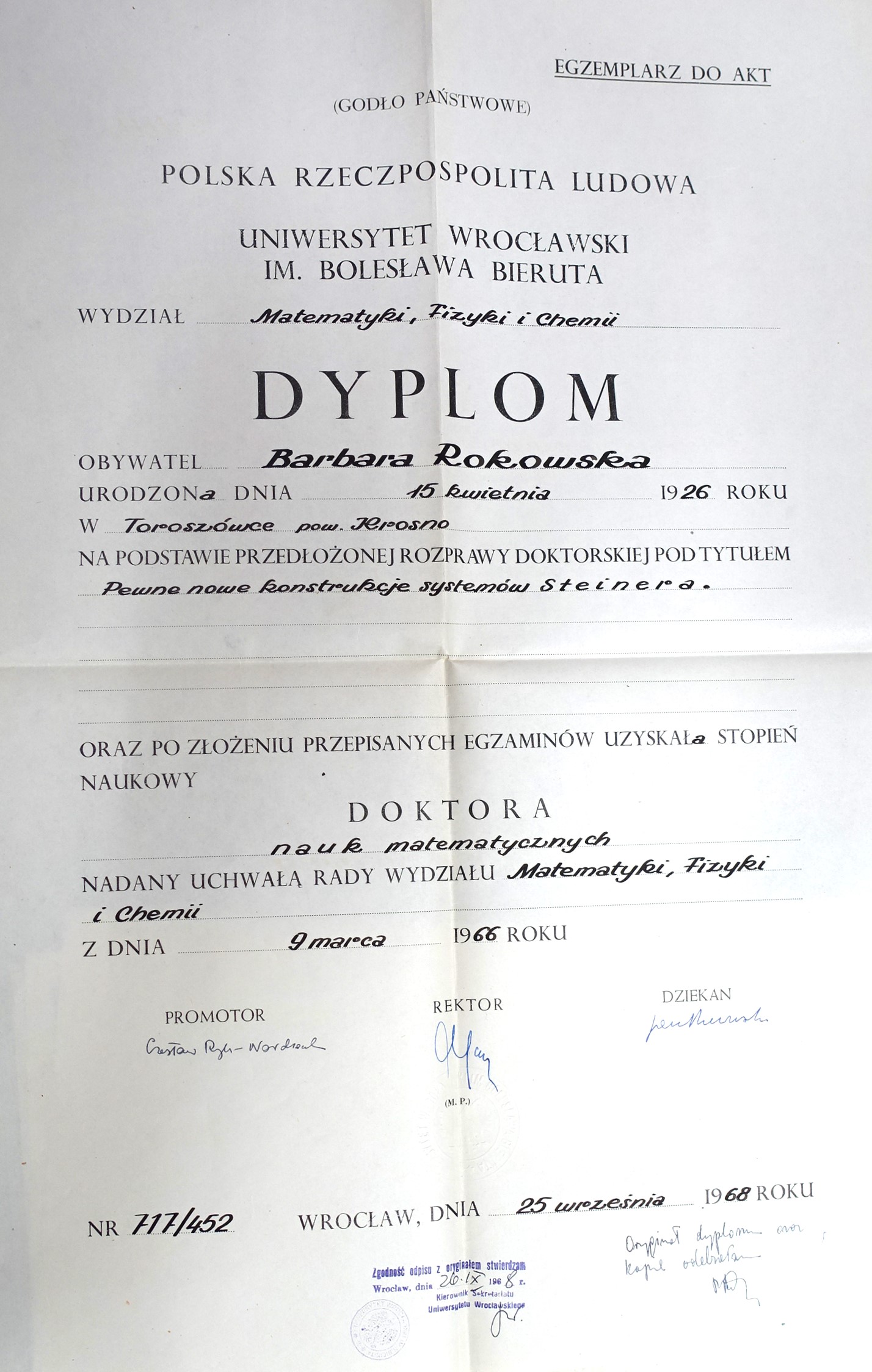}
       \caption{\small\label{BR_PhD2}PhD certificate.}
    \end{subfigure}
    \hfill
    \begin{subfigure}[b]{0.492\textwidth}
        \includegraphics[width=0.99\textwidth]{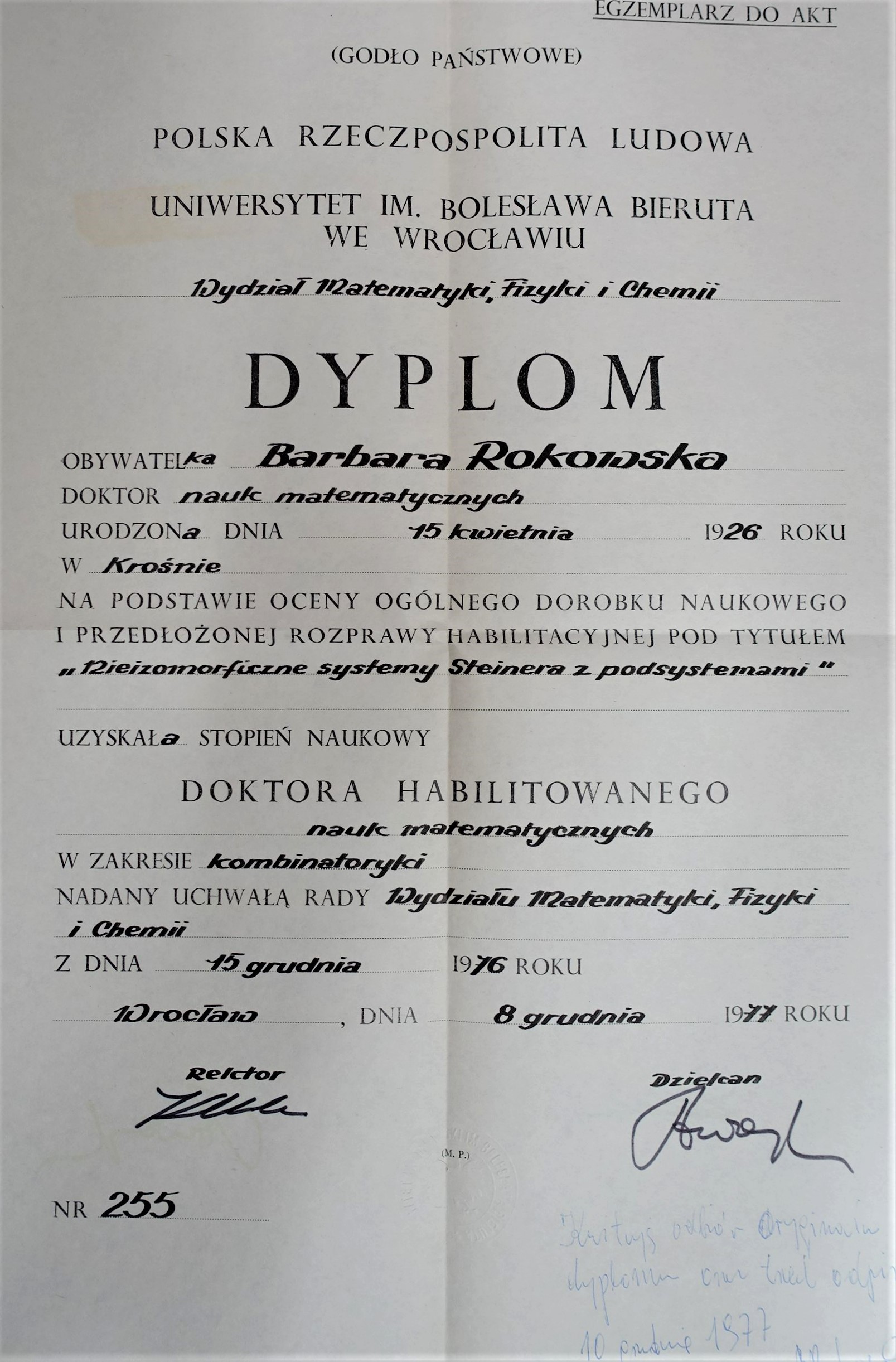}
      \caption{\small\label{BR_DyplHab} Habilitation certification.}
    \end{subfigure}
    \caption{\label{BR_Dypl} PhD and the post PhD degree (habilitation) certificates(v. \cite[mks Syg. Akt 34/500]{BRok1966:PhD}).}
~ 
\end{figure}

\section{\label{BRHab}Habilitation thesis topics.} The good reception of the results contained in the dissertation probably caused Rokowska to develop this line of research. We will devote some space here to discussing the fundamental purpose of creating different constructions of block systems. A new structure is interesting when its use leads to a new structure. A structure with the same parameters as a known structure is new if each renumbering of the elements of the set (each permutation) leads to an identical system of blocks. Each such permutation, therefore, is an automorphism in the sense that it is a \emph{permutation of the elements of the set that transforms the arrangement of Steiner triples\footnote{Steiner triple systems ({\STS}s)\label{KSzSTS}} into itself, preserving the structure of all triples}. It is known that on a set $V$ having $v$ elements, one can define $v!$ permutations. However, if in the considered set of all permutations we look for those that preserve the structure of the block system, we will notice that for most of the cardinality $v$ of the $V$ set, for which the block system exists, only some permutations preserve the structure of the block system. We will illustrate this for systems of Steiner triples.

\section{\label{sBRPhDSuper}PhD theses supervised by Dr. Hab. Barbara Rokowska.} The research topics dealt with by Rokowska were also carried out by her collaborators. Findings conducted under Rokowska's by PhD students are listed at Table~\ref{tab:TheDisertations}.

Barbara Rokowska's PhD students in chronological order are Krystyna \citeauthor{Wil1981:PhD}~(\citeyear{Wil1981:PhD}), Bogusław \citeauthor{Paw1984:PhD}~(\citeyear{Paw1984:PhD}), Jolanta \citeauthor{Arszynska1987:PhD}~(\citeyear{Arszynska1987:PhD}),  Mirosław \citeauthor{Bak1991:PhD}~(\citeyear{Bak1991:PhD}), Dariusz \citeauthor{Piegatowski1993:PhD}~(\citeyear{Piegatowski1993:PhD}),  (v.~Table~\ref{tab:TheDisertations} for further details). All doctoral dissertations (doctoral procedures, proceedings for granting the degree of doctor of mathematics) were conducted by the Scientific Council of the Institute of Mathematics of the Wrocław University of Technology.

\begin{table}[h!tb]
	\centering\scriptsize 
		\begin{tabular}{|c|l|p{15em}|p{9em}|c|} \hline\hline
		No&{ Author}&Title of a thesis&Reviewers&Year \\ \hline\hline
1&Krystyna  \citeauthor{Wil1981:PhD}~&{\tiny Roz\-dziel\-ne sys\-temy czwó\-rek $\text{RB}(v,4,3)$}&\multirow{3}{7em}{\tiny Mieczysław~Borowiecki\\ Wiktor~Marek\\ Zbigniew~Romanowicz}&\citeyear{Wil1981:PhD}\\ &    && &\\ 
&    && &\\ \hline
2&Bogusław  \citeauthor{Paw1984:PhD}&\tiny Prawie rozdzielne systemy blokow. &\multirow{2}{7em}{\tiny Jerzy~Płonka\\ Zbigniew~Romanowicz} &\citeyear{Paw1984:PhD}\\ 
 &    &&&\\\hline
3&Jolanta  \citeauthor{Arszynska1987:PhD}&\tiny Pewne własności symetrycznych systemów i kodów liniowych.&\multirow{2}{7em}{\tiny Mieczysław~Borowiecki\\ Zbigniew~Romanowicz}&~\citeyear{Arszynska1987:PhD}\\ &    && &\\ \hline
4&Mirosław  \citeauthor{Bak1991:PhD}&\tiny O rozdzielnych i prawie rozdzielnych systemach bloków.& \multirow{2}{7em}{\tiny Mieczysław~Borowiecki\\ Jerzy Płonka}&\citeyear{Bak1991:PhD}\\ 
&    && &\\ \hline
5&Dariusz  \citeauthor{Piegatowski1993:PhD}&\tiny Łuki zupełne w systemach trójek Steinera.&\multirow{2}{7em}{\tiny Mieczysław~Borowiecki\\ Jerzy Płonka }&\citeyear{Piegatowski1993:PhD}\\ &    && &\\ \hline
\end{tabular}
\caption{\label{tab:TheDisertations}Doctoral theses supervised by Dr. Hab. Barbara Rokowska.}	
\end{table}\normalsize

\section{Summary}
The described history of developing mathematical interests and pursuing a professional career is not typical. The subject, which is now very popular and important for applications, was treated by mathematicians in the past as a training ground for logical thinking. It may be worth revising this approach, because interesting mathematical problems, logical puzzles, usually inspire. Sometimes to use them. 

\appendix
\section{Supplementary material.}
\subsection{\label{KSzBRBio}Barbara Rokowska (1926--2012).} 
\begin{figure}[th!]
\centerline{\includegraphics[width=0.8\textwidth]{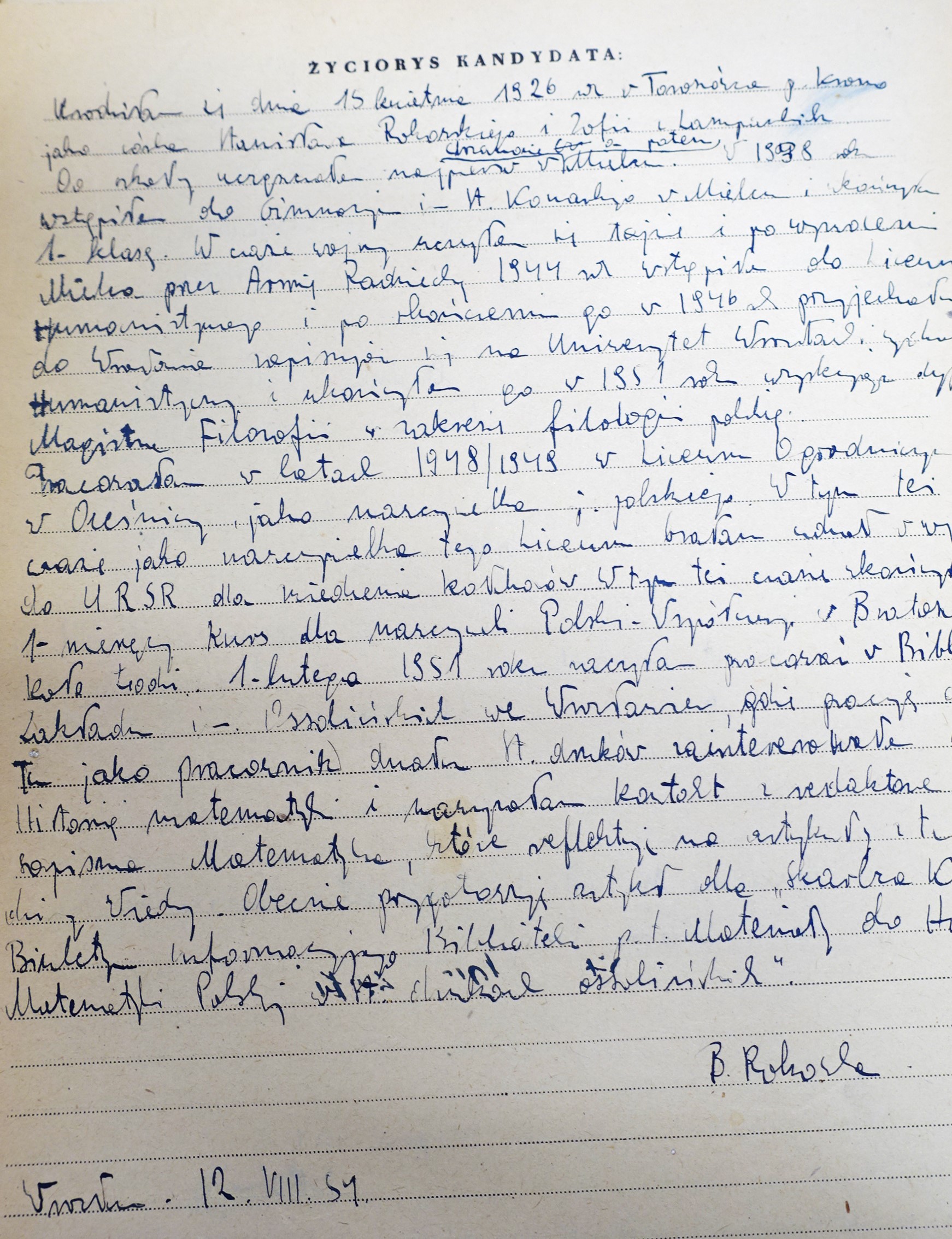}}
\caption{\label{BRCV1}Handwritten biography dated 12.viii.1954. (\cite[rks poz. 3]{BRokowska1954:UWr}).}
\end{figure}

She was born on April 15, 1926 in a small village, Toroszówka, near Krosno on Wisłoka\footnote{The most important events and facts of Rokowska's life were recorded by her, some more extensively described in a biographical book (v.~\citeauthor{Rok1996:Jezus}~(\citeyear{Rok1996:Jezus}). Those included in this article are also corroborated in the archival documents of the \UWr, \citeauthor{BRokowska1959:PWr}~(\citeyear{BRokowska1959:PWr}) and \citeauthor{BarRok2012:Ludzie}~(\citeyear{BarRok2012:Ludzie}).},  as the second child of Stanisław and Zofia (née ~Lamparska) Rokowski. Her father was an educated gardener from Wadowice, and took part in World War I in Pilsudski's legions. After the war, he remained in the army for a while, stationed in Nieszawa where he met his future wife, Zofia Lamparska (v. Fig~\ref{BRCV1}-- rks pos. 3 of \citeauthor{BRokowska1954:UWr}~(\citeyear{BRokowska1954:UWr}), p. 8 of \citeauthor{Rok1996:Jezus}~(\citeyear{Rok1996:Jezus})). 
\begin{quotation}\textsl\small I was not yet six years old when I was sent to a rural two-year school.
\end{quotation}(v. p. 9 of \citeauthor{Rok1996:Jezus}~(\citeyear{Rok1996:Jezus})). 

I studied at this school for one year and then at a school (grades 2-4) in Chrzanów until 1935. That was also the year when the Rokowski family moved to Mielec. She lived there longer. By World War II, she had completed the first grade of junior high school. When the war broke out, she was 13 years old. The beginning of the German aggression against Poland finds her and her mother outside Mielec, in her mother's family home in Nieszawa. After the end of the front operations, they return to Mielec. After some time, her father and brother were found, but throughout the war she was only with her mother. Her father participated in the resistance movement. At first as a courier for escaping soldiers to France. This activity of her father almost did not end with his arrest. 
\begin{wrapfigure}{l}{0.35\textwidth}
  \vspace{-8mm}
  \begin{center}
    \includegraphics[width=0.34\textwidth]{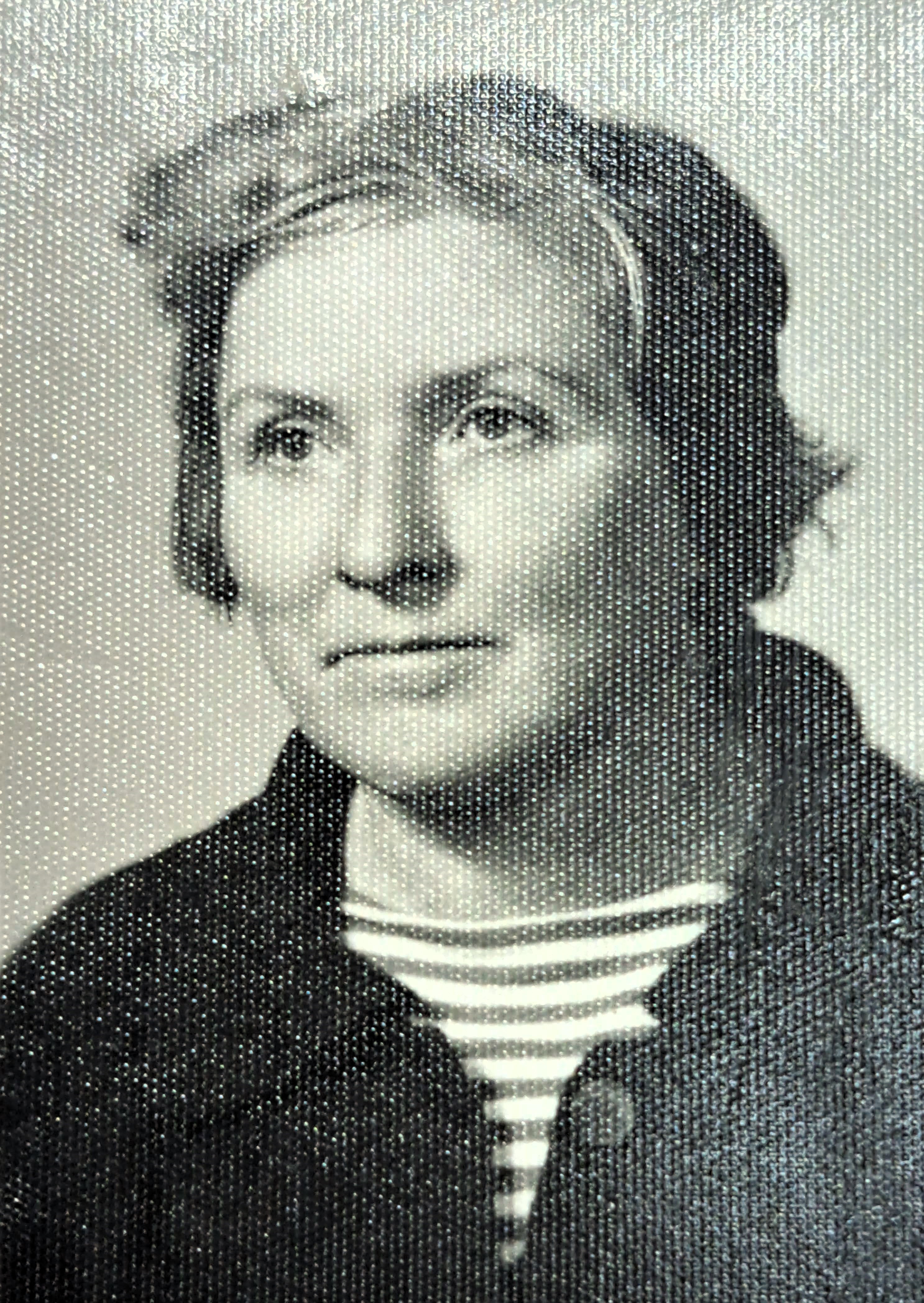}
		\includegraphics[width=0.335\textwidth]{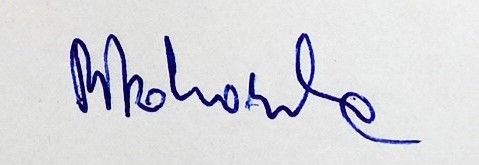}
   \caption{\small\label{KSzBRSign}Barbara Rokowska}
  \end{center}
  \vspace{-40pt}
\end{wrapfigure} 
In any case, he could not stay with his family because he was wanted by the Gestapo. Within the framework of the so-called sets, she worked on the material of 3 classes, from the 2nd to the 4th grade. Secret education was well organized in the region. Every year an examination committee came to the learners. They confirmed the results of the exam with a certificate, on which the pupil's data were hidden under a pseudonym. In 1944, after the Soviets captured Mielec and the Germans fled across the Vistula River, a regular Polish school was established. Although after the war her parents decided to move to the territories allotted to Poland ("the recovered lands"), they left their daughter in Mielec, thinking that she should not change schools until she finished school and took her matriculation exam. Thanks to the help of kind people, she survived the difficult winter of 1945/46 and passed her matriculation exam. After graduation, she went to her parents in Świdnica.

In 1946, at the age of 20, she began studying Polish philology at the University of Wrocław, which she graduated in 1951. During her studies, she is gainfully employed at a horticultural high school in Olesnica as a Polish language teacher (1948-49). Since 1950, she worked at Ossolineum Library in various departments. Among other things, in 1953 she was a member of the team preparing the catalog of Polonica of the 16th century (under the Kazimierz Zathey's direction) (v. \citeauthor{Tysz2000:Ossoliniana}~(\citeyear{Tysz2000:Ossoliniana})). 

\subsection{\label{BRMS}Mathematical Studies.} In 1954, influenced by friends, against her parents' suggestions and without their support, especially financially, she decided to study mathematics as a full-time student. This meant giving up her full-time job at the time. She originally explained her decision by her interest in the history of mathematics. Her early publications in the history of mathematics are BR(\ref{Rok1955:His1})H--BR(\ref{Rok1962:His9})H. After a while, her results made her eligible for a scholarship. She also worked as a text editor at the State Mathematical Institute in Wrocław. Her thesis, ``On the estimation of the integral $\int_0^{2\pi}\mid\cos(n_1 x)+\cos(n_2 x)+\ldots+\cos(n_k x)\mid dx$''  was related to the editorial job. It was described by Stanisław Hartman in his review of the thesis.
\begin{quote}[Doc.~Dr. Stanisław Hartman] Prof. Zahorski posed the question of the estimation of this integral depending on $n_k$ (with an increasing sequence $n_1<n_2<\ldots n_k)$ and expressed the conjecture that it is $\textbf{O}(\log n_k)$.

This conjecture was rebutted by P.~Erd\H{o}s, giving two pieces of evidence that the $\textbf{O}(\sqrt{n_k})$ estimate is overall better. He sent these proofs to the editors of the journal Colloquium Mathematicum in sketched form and left it to the editors to organize them into a printable manuscript. This task was given to student B.~Rokowska as the topic of her master's thesis. Since the two Erd\H{o}s  proofs operate with two completely different methods, each of which requires peculiar knowledge and orientation, the author\footnote{of master's thesis Barbara Rokowska} had to familiarize herself with several works. In addition to producing a compact manuscript to fill in the gaps and insinuations left by Erd\H{o}s.
\end{quote} 
\begin{figure}[th!]
    \centering
    \begin{subfigure}[b]{0.494\textwidth}
    \includegraphics[width=\textwidth]{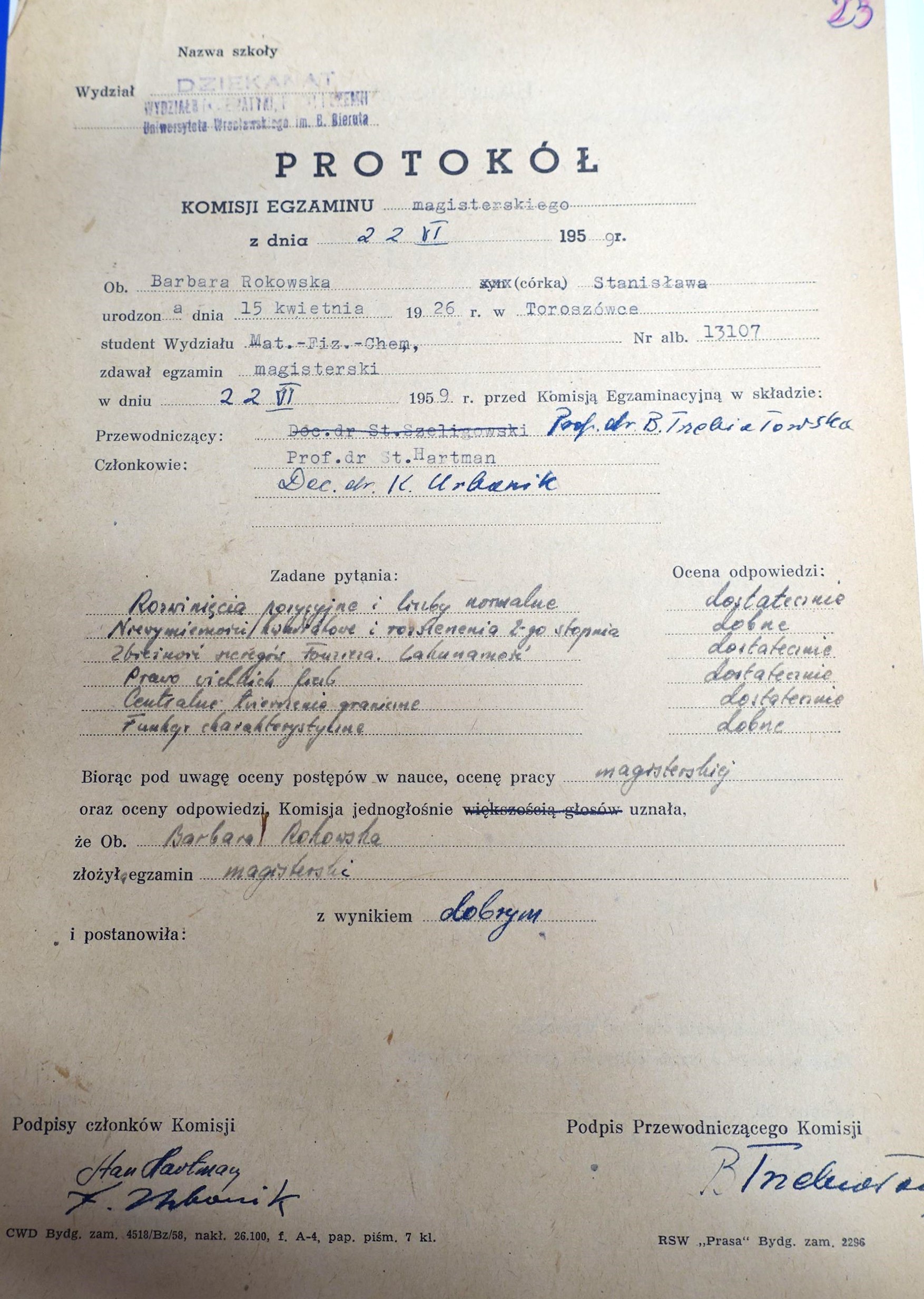}
       \caption{\small\label{BRProtoMgr}  The protocol of the exams.}
    \end{subfigure}
    \hfill
    \begin{subfigure}[b]{0.494\textwidth}
        \includegraphics[width=0.962\textwidth]{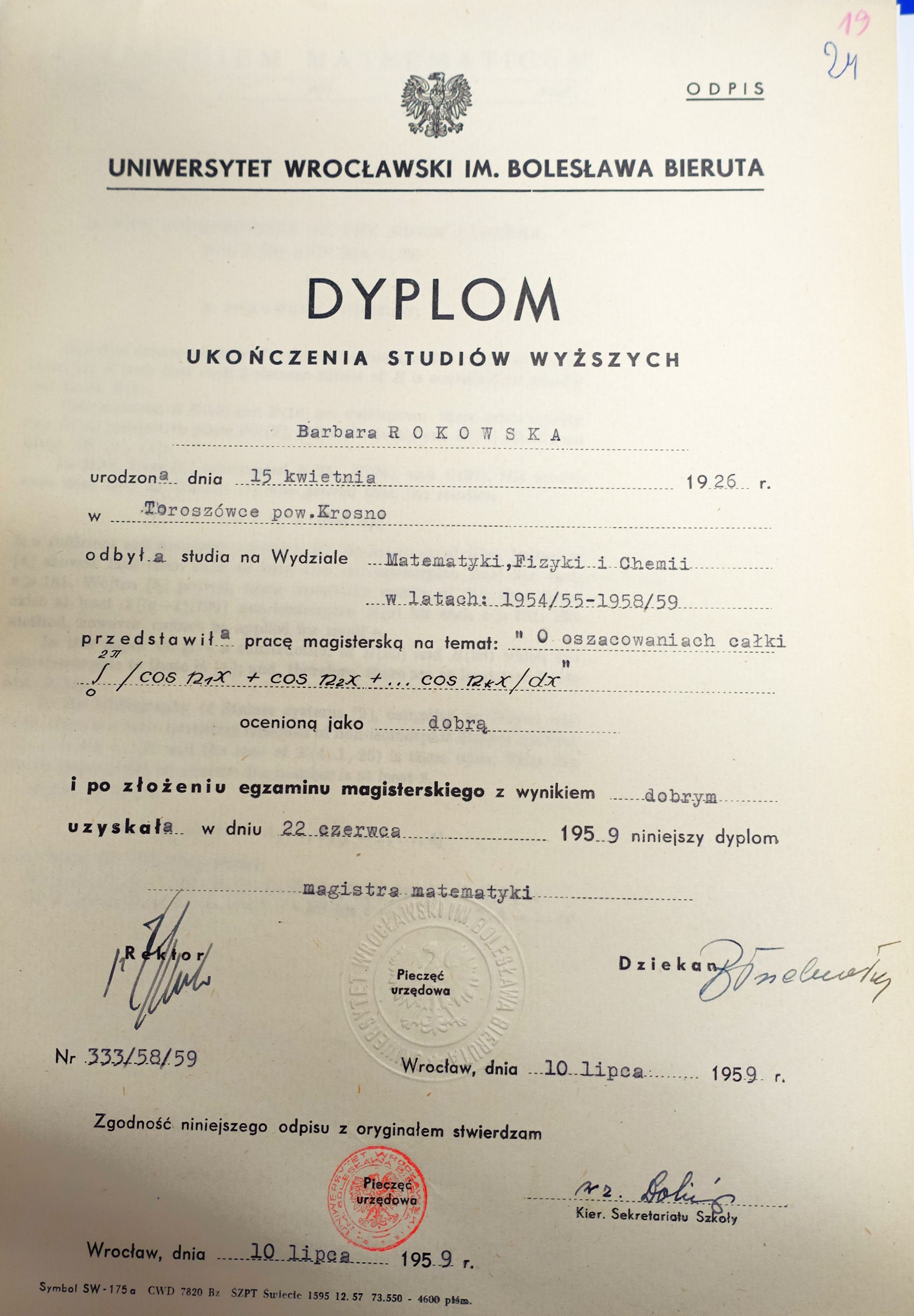}
      \caption{\small\label{BRDyplMgr}The master of science diploma.}
    \end{subfigure}
    \caption{\label{BRMgr} Master thesis records (v. \cite[mks 23,24 AUWr]{BRokowska1954:UWr})}
~ 
\end{figure}
\begin{figure}[th!]
    \centering
    \begin{subfigure}[b]{0.494\textwidth}
    \includegraphics[width=\textwidth]{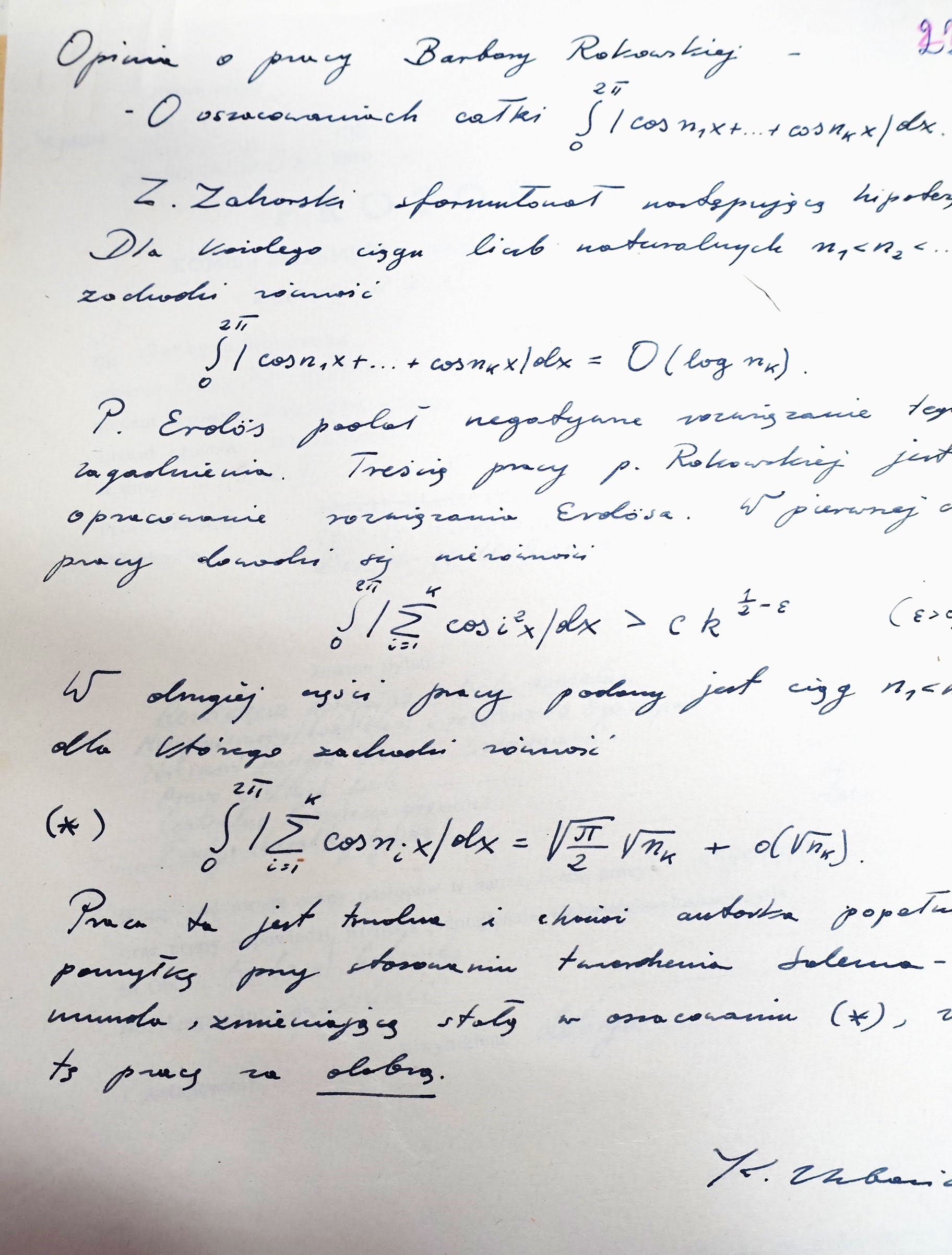}
       \caption{\small\label{BRRevKUMgr} Urbanik's review of the thesis.}
    \end{subfigure}
    \hfill
    \begin{subfigure}[b]{0.494\textwidth}
        \includegraphics[width=0.962\textwidth,height=0.4\textheight]{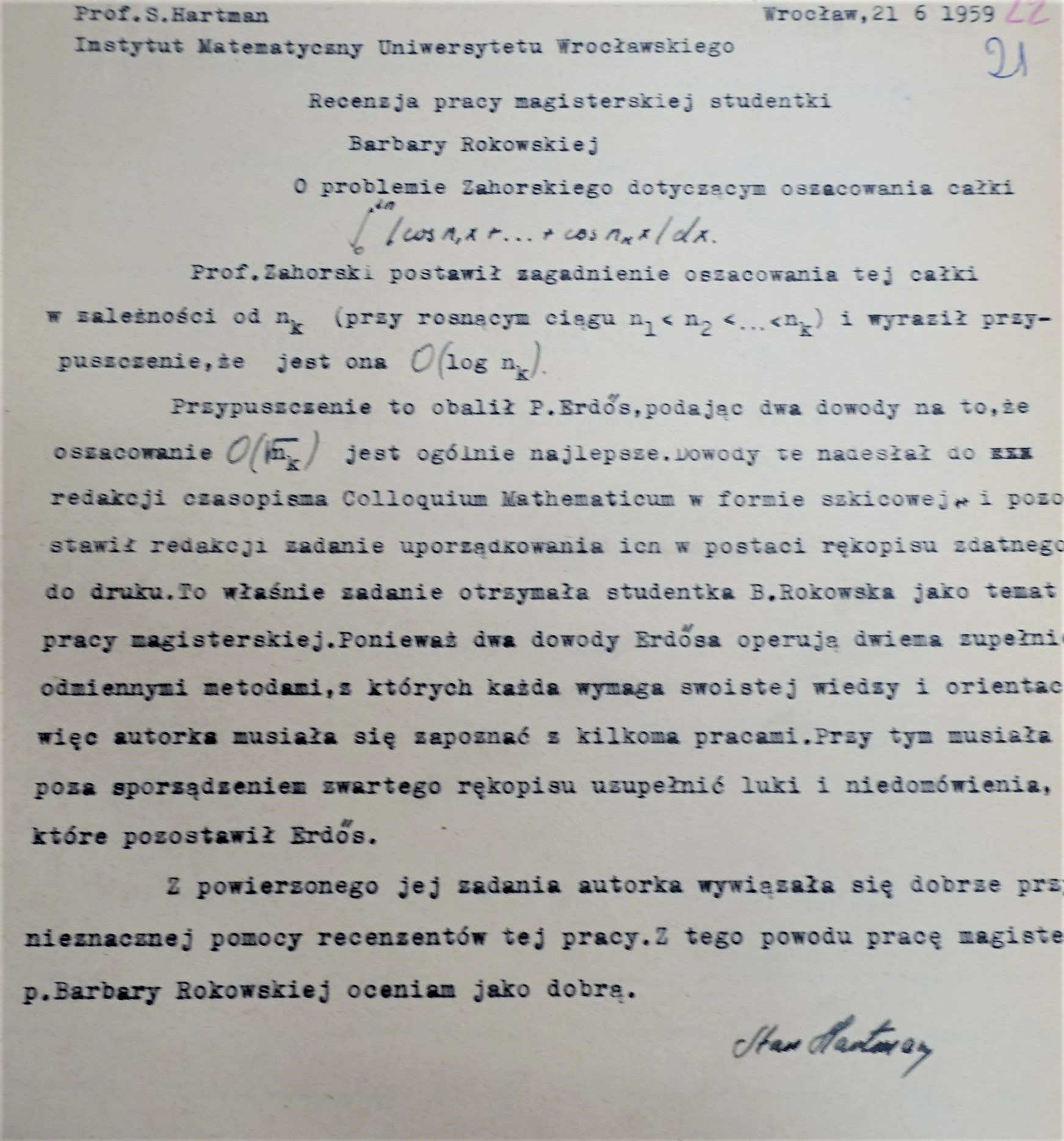}
      \caption{\small\label{BRRevSHMgr} Hartman's review of the thesis.}
    \end{subfigure}
    \caption{\label{BR2RevMgr} Master thesis records (v. \cite[mks 23,24 AUWr]{BRokowska1954:UWr})}
~ 
\end{figure}
\begin{figure}[th!]
\centerline{\includegraphics[width=9cm]{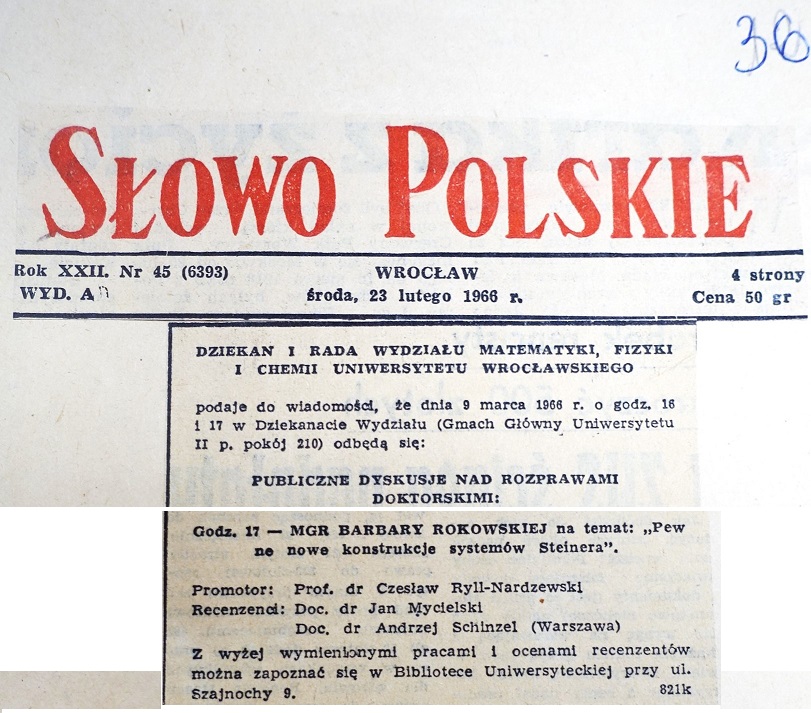}}
\caption{\label{BRPhD3}Ogłoszenie w \emph{Słowie Polskim} o terminie obrony doktoratu (v. \cite[mks poz. 36]{BRok1966:PhD})}
\end{figure}

\subsection{\label{BRlabor}Her first employment as a mathematician} was at the Wroclaw University of Technology\footnote{Politechnika Wrocławska (\PWr)\label{BRPWr}} as an assistant in the Chair of Mathematics. She began work on October 1, 1959 with the start of the 1959/60 academic year (\cite[poz.~4]{BRokowska1959:PWr}). As it later turned out, she remained there throughout her further professional life. Her work at the Polytechnic mainly boiled down to teaching courses in mathematics, lectures and exercises in mathematical analysis and algebra to students of various technical majors. Longer employment was possible when, in addition to teaching, the employee conducted research and published the results. Rokowska pursued her passion for research by continuing to study block systems (see definitions in Section~\ref{sBRPhD}).  She already had experience in editing mathematical papers at that time, but she did not publish the results constituting her dissertation before starting the doctoral procedure. The territorially closest Scientific Council with the authority to conduct a doctoral procedure for Rokowska was at her Alma Mater. She also submitted her dissertation there. Her supervisor was Czesław Ryll-Nardzewski. It is noteworthy that the topic of her research was not related to the supervisor's main research interests. 

It is worth mentioning here that she was a highly regarded lecturer. When the Faculty of Fundamental Problems of Technology was established in 1968, the first students of the Faculty were taught calculus by Rokowska, among others. Organizationally, the Faculty consisted of two parts: Institute of Mathematics and Theoretical Physics (\IMiFT\label{BR_IMiFT}) and the second one, Institute of Physics. Dr. Rokowska was in the first one.     

She was an active member of Polish Mathematical Society\footnote{V. record in \citeauthor{WM2316:PTM}~(\citeyear{WM2316:PTM}):Rokowska Barbara, mgr, Wrocław, ul. Grabiszyńska 29/1.}.

\subsection{\label{BRPhDproc}Doktor of Philosophy Dissertation} entitled. “\emph{Several new constructions of Steiner systems}” was defended on March 9, 1966 at the Faculty of Mathematics, Physics and Chemistry of \UWr. The doctoral thesis was promoted by Czeslaw Ryll-Nardzewski (2026-2015)\footnote{V.~\citeauthor{Kwapien2015:CRN}~(\citeyear{Kwapien2015:CRN})}. The reviewers were Professors Jan Mycielski\footnote{Jan Mycielski is a Polish-American mathematician, logician and philosopher, a professor emeritus of mathematics at the University of Colorado at Boulder. He is known for contributions to graph theory, combinatorics, set theory, topology and the philosophy of mathematics.} and Andrzej Schinzel(1904-1969). Quoting the conclusions of the reviews, the author worked on old issues, whose importance has recently been rediscovered. An important part of the dissertation [BR\ref{Rok1966:TPhD}] was published by Rokowska in \CM (v.~[BR\ref{Rok1967:PhD1}]).

\subsection{Rokowska's habilitation procedure. } On December 9, 1974, the Scientific Council of the Faculty of Mathematics, Phisics and Chemistry \UWr\  appointed a committee for the \emph{postdoctoral dissertation} (habilitation) of Dr. Barbara Rokowska. The proposal was referred by Prof. Władysław Narkiewicz. The commission was appointed, consisting of professors: Czesław Ryll-Nardzewski, Ludwik Borkowski and Władysław Narkiewicz. On the same session of the Scientific Council, the commission proposed as reviewers professors: Jerzy Browkin from Warsaw University, Stanisław Hartman from \IM\label{BR_IM} \PAN\label{BR_PAN}, and Leon Jeśmianowicz from Mikołaj Kopernik University in Toruń. 

\paragraph{Reviews} of the dissertation were received with a long delay. The reasons for this delay are explained in his post (the review) by Professor Hartman. In an effort to find out whether the content of the dissertation is correct, he put a lot of time into studying the first part of the dissertation together with Dr. Rokowska (v.~App.~\ref{BRhabRecSH}). At the time, this part was submitted for publication in \emph{Colloquium Mathematicum}, where Hartman was one of the editors. Although the reviewer asked by the Editorial Board recommended the work for publication, however Hartman's inquisitiveness caused the work to be significantly rewritten, which, according to the reviewer, greatly improved its message. In his review, Hartman did not make a motion to grant Dr. Rokowska a postdoctoral degree. He stressed that the results of the presented dissertation are from , “finistic mathematics,” which, while not exclusive, also the lack of conclusions (theorems, proofs, methods) in the presented dissertation, from the broader area of mathematics, does not allow him to support the request of the interested party.

Let us recall that the term "finitistic" refers to finitism, a position in the philosophy of mathematics that only accepts mathematical operations and objects that are finite or can be constructed in a finite number of steps. Finitistic approaches exclude the use of infinity as an actual mathematical object.

In particular, within finitism:
\begin{itemize}\itemsep2pt \parskip2pt \parsep0pt
\item Only finite and concrete proofs and constructions are accepted.
\item The use of infinity (e.g., infinite sequences) is considered problematic or unacceptable.
\end{itemize}
 Finitism was notably developed by mathematicians like David Hilbert, who sought to establish formal foundations for mathematics without relying on infinity.

It was late 1976 when the three reviews went to committee. In the meantime, a legal flaw was raised at the Scientific Council meeting. One of its members had changed his place of employment, causing him to cease to be a member of the board considering the request. On September 23, 1976, the composition of commission was supplemented by professor Józef Łukaszewicz, due to the change in Borkowski's affiliation. After receiving all the reviews, the Committee (all 4 members) held a meeting on October 30, 1976 and found that two of them were positive, while one, the review by Prof. Hartman, was negative. The Committee decided to admit Dr. Barbara Rokowska to the habilitation colloquium.
\begin{quote}[Excerpt from “Application for approval of habilitation...” (v.~\cite[Card 21, 18]{Rokowska1976})
On October 30, 1976 the committee analyzed the reviews and made a motion to proceed further, that is, to hold a colloquium before the Scientific Council . 
\end{quote}
\paragraph{The colloquium} was on December 15, 1976. Prof. Wajda presided and doc Krzywicki took the minutes. Out of 52 members of the Scientific Council entitled to vote and 36 present, and 3 reviewers, 35 people voted in favor of the colloquium and 4 abstained, while 34 people voted in favor of awarding the postdoctoral degree and 5 abstained.  

\subsection{Professional and social activity.} Barbara Rokowska chose her subject well. She was able to interest others with her research, including young mathematicians looking for challenges leading to a doctorate. The first doctoral student to be promoted was Krystyna Wilczyńska (defense in 1981). Earlier, in 1974, Krystyna Pukanow\footnote{Without habilitation, Dr. Rokowska could not act as a supervisor or assistant supervisor in accordance with the law in force at that time.} defended her doctoral dissertation. The theses of her doctorate partly concerned topics formulated by Rokowska, and also analyzed together with her. Pukanow's doctoral supervisor was Doc. Dr. Lucjan Szamkołowicz. Since habilitation, she has promoted five doctors. The last of them was Dariusz Piegatowski (defense in 1993).

At the Wrocław University of Technology, she went through all stages of an academic teacher's career. Employed in 1959 as an assistant, in 1961 she was promoted to senior assistant. She received the position of assistant professor in 1968. From 1977 she had the position of associate professor\footnote{The Minister signed the appointment on 23 August 1977 for the position of associate professor effective 1 September 1977}, and in 1990 she received the position of associate professor at the Wrocław University of Science and Technology. She was granted a pension on 26 June 1992, but remained in full-time employment until 30 September 1996. In the following years she still taught part-time until 2002.

From her autobiographical books (v.~\cite{Rok1996:Jezus,Rok2005:Jezus}) and archives of the Wrocław University of Science and Technology we know that she had a son, Jacek, who was born in 1961. His upbringing was a significant burden and concern, which caused that from a certain point in her life she helped people who were lost and needed friendly support. She devoted herself to social activities, worked, among others, in the Crusade for the Liberation of Man, and ran a day care center for children from difficult homes at the parish of St. Michael in Wrocław.

\paragraph{Ph.D. Barbara Rokowska died on June 1, 2012 after a long illness. } The funeral mass was held on Wednesday, June 6, 2012, Wednesday 10:40, in the Church of St. Michael the Archangel at ul. Bolesława Prusa. She was buried at the Osobowice Cemetery [field 134, row 6, grave 192]

\section{Barbara Rokowska's research publications.} 
  
\subsection{History of mathematics}
{\small
\begin{enumerate}[{\hspace{-2em}BR(1)H}]\itemsep1pt \parskip2pt \parsep1pt
\providecommand{\natexlab}[1]{#1}
\providecommand{\url}[1]{\texttt{#1}}
\expandafter\ifx\csname urlstyle\endcsname\relax
  \providecommand{\doi}[1]{doi: #1}\else
  \providecommand{\doi}{doi: \begingroup \urlstyle{rm}\Url}\fi

\item\label{Rok1955:His1}
B.~Rokowska.
\newblock Dickstein i historycy.
\newblock \emph{Matematyka}, 8\penalty0 (1/2):\penalty0 5--7., 1955.

\item\label{Rok1956:His2}
B.~Rokowska.
\newblock Podręczniki matematyki z okresu {K}omisji {E}dukacji {N}arodowej.
\newblock \emph{Matematyka}, 9\penalty0 (2):\penalty0 20--27., 1956.

\item\label{Rok1957:His3}
B.~Rokowska.
\newblock Abraham {S}tern z {H}rubieszowa, pierwszy polski konstruktor maszyn
  rachunkowych.
\newblock \emph{Matematyka}, 10\penalty0 (3):\penalty0 13--17., 1957.

\item\label{Rok1958:His4}
B.~Rokowska.
\newblock {Józef Puzyna}.
\newblock \emph{Matematyka}, 9\penalty0 (4-6):\penalty0 4--7., 1958.

\item\label{Rok1960:His5}
B.~Rokowska.
\newblock Niewęgłowscy. Ojciec i syn.
\newblock \emph{Matematyka}, 13\penalty0 (4):\penalty0 140--144., 1960{\natexlab{b}}.

\item\label{Rok1960:His6}
B.~Rokowska.
\newblock O liczbach postaci $2^n3^m + 1$.
\newblock \emph{Matematyka}, 13\penalty0 (5):\penalty0 253--258., 1960{\natexlab{c}}.

\item\label{Rok1960:His7}
B.~Rokowska.
\newblock Stanisław Zaremba.
\newblock \emph{Matematyka}, 13\penalty0 (5):\penalty0 328--332., 1960{\natexlab{d}}.

\item\label{Rok1961:His8}
B.~Rokowska.
\newblock O kwadratach łacińskich.
\newblock \emph{Matematyka}, 14\penalty0 (2):\penalty0 69--70, 1961.

\item\label{Rok1962:His9}
B.~Rokowska.
\newblock Władysław Folkierski.
\newblock \emph{Matematyka}, 15\penalty0 (5):\penalty0 263--266., 1962.

\end{enumerate}
}

\subsection{Research papers}
{\small
\begin{enumerate}[{\hspace{-2em}BR1}]\itemsep1pt \parskip2pt \parsep1pt
\providecommand{\natexlab}[1]{#1}
\providecommand{\url}[1]{\texttt{#1}}
\expandafter\ifx\csname urlstyle\endcsname\relax
  \providecommand{\doi}[1]{doi: #1}\else
  \providecommand{\doi}{doi: \begingroup \urlstyle{rm}\Url}\fi

\item\label{Rok1959:First}
B.~Rokowska.
\newblock On periodic sequences of natural numbers.
\newblock \emph{Wiadom. Mat. (2)}, 3:\penalty0 41--43 (1959),
  1959{\natexlab{a}}.
\newblock ISSN 0373-8302.
\newblock \doi{10.14708/wm.v3i1.2679}.
\newblock \MR{115961}, \ZBL{0095.26104}.

\item\label{RokSchi1960:Erdos}
B.~Rokowska and A.~Schinzel.
\newblock Sur un probl\`eme de {M}. {E}rd\H{o}s.
\newblock \emph{Elem. Math.}, 15\penalty0 (4):\penalty0 84--85,
  1960.
\newblock ISSN 0013-6018.
\newblock URL
  \url{http://resolver.sub.uni-goettingen.de/purl?PPN378850199_0015}.
\newblock \MR{117188}; \ZBL{0089.26603}.

\item\label{Rok1966:TPhD}
B.~Rokowska.
\newblock \emph{{Pewne nowe konstrukcje systemów Steinera}}.
\newblock Rozprawa doktorska, Wroclaw University. Faculty of Mathematics,
  Physics and Chemistry, Wroc{\l}aw, 1966.

\item\label{Rok1967:PhD1}
B.~Rokowska.
\newblock Some new constructions of {$4$}-tuple systems.
\newblock \emph{Colloq. Math.}, 17:\penalty0 111--121, 1967.
\newblock ISSN 0010-1354.
\newblock \doi{10.4064/cm-17-1-111-121}.
\newblock \MR{0213253}; \ZBL{0158.01404}

\item\label{Rok1968:RokSzm}
B.~Rokowska, Z. Szmorliński.
\newblock Obwody powierzchniowe maszyn elektrycznych a przekształcenie trójkąt-gwiazda.
\newblock \emph{Archiwum Elektroniki}, 17:\penalty0 41--46, 1968.

\item\label{Rok1971:PhD2}
B.~Rokowska.
\newblock Some remarks on the number of different triple systems of {S}teiner.
\newblock \emph{Colloq. Math.}, 22:\penalty0 317--323, 1971.
\newblock ISSN 0010-1354.
\newblock \doi{10.4064/cm-22-2-317-323}.
\newblock \MR{286681}; \ZBL{0218.05008}.

\item
B.~Rokowska.
\newblock On the number of different triple systems of {Steiner}.
\newblock Prace Nauk. {Inst}. {Mat}. {Fiz}. Teor., {Politechniki} {Wrocław}.,
  {Ser}. {Studia} {Materialy} 1972, {No}. 6, 41-57 (1972)., 1972{\natexlab{a}}.
\newblock \MR{354407}; \ZBL{0254.05005}.

\item
B.~Rokowska.
\newblock Some new configurations {$C(k,l,n,\lambda)$}.
\newblock Prace {Nauk}. {Inst}. {Mat}. {Fiz}. {Teor}., {Politechniki}
  {Wroclaw}., {Ser}. {Studia} {Materialy} 4, 61-64 (1972)., 1972{\natexlab{b}}.
\newblock ISSN 0324-9581.
\newblock \MR{344142}; \ZBL{0362.05050}.

\item\label{Rok1973:STS73}
B.~Rokowska.
\newblock On the number of non-isomorphic {Steiner} triples.
\newblock \emph{Colloq. Math.}, 27:\penalty0 149--160, 1973.
\newblock ISSN 0010-1354.
\newblock \doi{10.4064/cm-27-1-149-160}.
\newblock \MR{325412}; \ZBL{0254.05006}.

\item
B.~Rokowska.
\newblock \emph{{Nieizomorficzne systemy trójek Steinera}}.
\newblock Rozprawa habilitacyjna, Wroclaw University. Faculty of Mathematics,
  Physics and Chemistry, Wrocław, 1976{\natexlab{a}}.

\item
B.~Rokowska.
\newblock On resolvable {BIBD}.
\newblock Graphs, {Hypergraphs}, {Block} {Syst}.; {Proc.} {Symp.} Comb. {Anal}., {Zielona} {Góra} 1976, 221-225 (1976)., 1976{\natexlab{b}}.
\newblock \ZBL{0357.05015}.

\item
B.~Rokowska.
\newblock Non-isomorphic {Steiner} triples with subsystems.
\newblock \emph{Colloq. Math.}, 38:\penalty0 153--164, 1977{\natexlab{a}}.
\newblock ISSN 0010-1354.
\newblock \doi{10.4064/cm-38-1-153-164}.
\newblock \MR{0460146}; \ZBL{0377.05004}.

\item
B.~Rokowska.
\newblock A new construction of the block systems {$B(4,1,25)$} and {$B(4,1,28)$}.
\newblock \emph{Colloq. Math.}, 38\penalty0 (1):\penalty0 165--167, 1977{\natexlab{b}}.
\newblock ISSN 0010-1354.
\newblock \doi{10.4064/cm-38-1-165-167}.
\newblock \MR{0460147}; \ZBL{0377.05005}.

\item
B.~Rokowska and M.~Wojtas.
\newblock Nonisomorphic balanced incomplete block designs {$B(v, 1, 5)$}.
\newblock \emph{Prace Nauk. Inst. Mat. Politech. Wroc\l aw. Ser. Stud.
  Materia\l y}, 17\penalty0 (13):\penalty0 19--24, 1977{\natexlab{a}}.
\newblock ISSN 0324-9808.
\newblock \MR{505626}; \ZBL{0378.05012}.

\item
B.~Rokowska.
\newblock On resolvable designs.
\newblock \emph{Scientific Papers of the Institute of Mathematics of Wrocław
  Technical University. Studies and Research}, 19\penalty0 (14 Analiza
  Dyskretna):\penalty0 3--9, 1980.
\newblock ISSN 0324-9808.
\newblock Politechnika Wrocławska. Instytut Matematyki. Prace Naukowe. Seria
  Studia i Materiały. \ZBL{03724457}, \MR{586887}.

\item
B.~Rokowska.
\newblock On resolvable quadruple systems.
\newblock \emph{Discuss. Math.}, 4:\penalty0 43--50, 1981.
\newblock ISSN 0137-9747.
\newblock \ZBL{0511.05013}.

\item
B.~Rokowska.
\newblock On resolvable 6-tuples.
\newblock \emph{Discuss. Math.}, 6:\penalty0 15--17, 1983.
\newblock ISSN 0137-9747.
\newblock \ZBL{0551.05021}.

\item
B.~Rokowska and M.~Wojtas.
\newblock On resolvable designs.
\newblock \emph{Discuss. Math.}, 6:\penalty0 19--26, 1983.
\newblock ISSN 0137-9747.
\newblock \ZBL{0548.05010}.

\item\label{Rok1984:STS84}
B.~Rokowska.
\newblock Resolvable systems of {$8$}-tuples.
\newblock \emph{J. Statist. Plann. Inference}, 9\penalty0 (1):\penalty0
  131--141, 1984.
\newblock ISSN 0378-3758.
\newblock \doi{10.1016/0378-3758(84)90050-8}.
\newblock \MR{MR0735017}; \ZBL{0532.05005}.

\item
B.~Rokowska.
\newblock The construction of resolvable block-systems.
\newblock \emph{Colloq. Math.}, 49\penalty0 (2):\penalty0 295--297,
  1985{\natexlab{a}}.
\newblock ISSN 0010-1354.
\newblock \doi{10.4064/cm-49-2-295-297}.
\newblock \MR{830811}; \ZBL{0571.05007}.

\item
B.~Rokowska.
\newblock Construction of a some class of nearly resolvable block-systems.
\newblock \emph{Discuss. Math.}, 7:\penalty0 17--23, 1985{\natexlab{b}}.
\newblock ISSN 0137-9747.
\newblock \MR{852834}; \ZBL{0594.05012}.

\item
B.~Rokowska and K.~Wilczyńska.
\newblock Certain constructions of {S}teiner systems.
\newblock \emph{Discuss. Math.}, 7:\penalty0 25--29, 1985.
\newblock ISSN 0137-9747.

\item
B.~Rokowska and K.~Wilczyńska.
\newblock A base for a system of seven-element blocks.
\newblock \emph{Discuss. Math.}, 8:\penalty0 7--12 (1987), 1986.
\newblock ISSN 0137-9747.
\newblock \MR{932080}; \ZBL{0683.05011}.

\item
B.~Rokowska and K.~Wilczyńska.
\newblock Decomposition of a complete graph into hexagons.
\newblock \emph{Discuss. Math.}, 9:\penalty0 45--54 (1989), 1988.
\newblock ISSN 0137-9747.
\newblock \MR{1042460}; \ZBL{0719.05051}.

\item
B.~Rokowska.
\newblock Decomposition of a complete graph into hexagons.
\newblock volume~45, pages 141--148 (1991). 1990.
\newblock Graphs, designs and combinatorial geometries (Catania, 1989).
\newblock  \MR{1181003}; \ZBL{0761.05074}.

\item
B.~Rokowska and K.~Wilczy\'{n}ska.
\newblock On number of intersections in a symmetrical system.
\newblock \emph{Discuss. Math.}, 11:\penalty0 5--16 (1992), 1991.
\newblock ISSN 0137-9747.
\newblock \MR{1178354}; \ZBL{0758.05005}.

\item
B.~Rokowska and K.~Wilczyńska.
\newblock On the existence of {$(v,4,2)$}--perfect {M}endelsohn designs.
\newblock \emph{Discuss. Math.}, 14:\penalty0 5--13, 1994.
\newblock ISSN 0137-9747.
\newblock \MR{1323968}; \ZBL{0824.05006}.
\end{enumerate}
}

\section{\label{BRhabRecSH}Recenzja prof. Stanisława Hartmana rozprawy habilitacyjnej dr. Barbary Rokowskiej}
The following is a full translation of the review of Dr. Barbara Rokowska's postdoctoral dissertation by Prof. Stanislaw Hartman. Its various excerpts show the broad cultural context of combinatorial research in Poland. In an environment dominated by mathematicians focused on selected areas of research, it is difficult to convince the community of the importance of research in other areas. \\[2ex]

\leftline{Prof. Stanislaw Hartman\hfill	Wrocław, 6.IX.1976}

\paragraph{ Evaluation of scientific achievements and habilitation thesis of Dr. Barbara Rokowska.} The entire scholarly output of Dr. Barbara Rokowska, with the exception of popular works on the history of mathematics and one joint paper${}^\text{[BR\ref{Rok1968:RokSzm}]}$, concerns Steiner systems and their generalizations. The original Steiner triples for a given finite set $A$ are such a system of three-element subsets of it that every pair of elements from $A$ is contained in one and only one of the intermeetings. The study of such systems, and in particular their existence dependent on the count $n$ of the set $A$, has proved important in many combinatorial problems. The concept itself can be naturally generalized: instead of triples, one can introduce subsets of cardinality $k$, instead of pairs subsets of fixed powers of $\ell<k$, and demand that each $\ell$-element subset of the set A be contained in exactly $\lambda$ $(\lambda >1)$ from among the sets forming the system sought, which thus depends on four parameters; $n$, $k$, $\ell$, $\lambda$ and is called a configuration $C(n,\ k,\ \ell,\ \lambda)$. Steiner's theory of triples and configurations (we will also say “Steiner systems” for shorter) developed for more than 150 years and dealt mainly with the existence of systems with given parameters and the number of those that are significantly different, that is nonisomorphic, i.e., not formed one from another by permutation of the set $A$.

These two problems are related, in that the existence proofs are constructive, and the constructions leave some arbitrariness that can be used for estimating from below the number of non-isomorphic systems. This number behaves very capriciously, and is difficult to study already for Steiner triples, e.g. For $n<15$ there are at most two non-isomorphic systems of triples, when $n=15$ there are exactly $80$. The increase in the number of non-isomorphic systems of triples with increasing cardinality of the set A is very fast (Aleksiejev estimates it in 1974 from below by $\exp{\big(\frac{1}{12} (n^2-\epsilon)\log(n)\big)})$\footnote{V.~\citeauthor{Alekseev1974:STS}~(\citeyear{Alekseev1974:STS}).}, but irregular.  The greatest progress in more recent research of this type is due to Hanani. Various mathematicians are currently working on these issues, e.g. Doyen\footnote{V.~\citeauthor{Doy1973:Steiner}~(\citeyear{Doy1973:Steiner,Doy1976:Steiner}) and \citeauthor{DoyRos1973:BibSteiner}~(\citeyear{DoyRos1973:BibSteiner,DoyRos1978:BibSteiner,DoyRos1980:BibSteiner})}, the aforementioned Aleksiejev, Robinson and the Author of the habilitation thesis discussed here, her colleagues and students in Wrocław.

Dr. Rokowska's postdoctoral dissertation consists of 3 parts: I Non-isomorphic Steiner triples with subsystems (in typescript accepted for publication in Colloquium Mathematicum in 1975), II On the number of non-isomorphic Steiner triples (Coll. Math 1975), III On the number of different triple systems of Steiner (Scientific Papers of the Institute of Mathematics, Wroclaw University of Technology 1972). 

In Part I, the author estimates from below the number of nonisomorphic systems of Steiner triples that contain subarrangements constructed on a set of a given cardinality s, obviously smaller than the cardinality n of the set A. It turns out that this number grows immeasurably with s if only $s \equiv 1 \text{ or } 3 (\mod 6)$, which is the classical necessary and sufficient condition for the existence of Steiner triples on the set of cardinality s. The main theorem of Paper I is a substantial generalization of the estimates obtained by the Author in the previous two papers, which constitute Parts III and II of the habilitation dissertation (using the Author's numbering here). There, it is about systems of triples containing subarrays constructed on sets of the highest possible cardinality. The proofs of theorems in all three parts of the dissertation are constructive and complicated. This is especially true of the thesis I, where the proof of the main theorem is based on some Hanani's methods, but required the author to do long and difficult reasoning, to get into the construction of a system of triples at different cardinality n and to handle and use configurations more complex than triples.  It seems to me that the term “virtuosity” will not be an exaggeration for this proof. Let's add that in some cases the author was able to find the exact number of non-isomorphic systems studied. 

The reviewer thinks it his duty to present a rather peculiar situation that arose around the main part of the reviewed dissertation, i.e., paper I. This work was accepted for publication in Colloquium Math with too much haste on the basis of a review by a mathematician selected by the editors, who, it turned out, did not pay enough attention to the way it was written and, above all, to the completeness of the proofs. When I later took a closer interest in the work , I pointed out its great shortcomings, and from then on I began to act in a dual role; as a reviewer of the postdoctoral dissertation and as deputy editor of Colloquium Math, who should take care of the proper form and comprehensibility of published articles. I had no doubt that the results of the work of Part I were fit for publication, and the reasoning showed such erudition and skill of the author in this field that one gained confidence in the correctness. However, checking the correctness seemed impossible in the face of unexpected reductions and clearly unjustified premises. The manuscript was extremely far from what is called “self-contained” in English, and could be understood at most after studying various works, in particular Hanani's, which the author only vaguely refers to. This is not a matter of leaving out the proofs of theorems proven explicitly elsewhere, which is of course allowed, but the lack of explanation of where certain important passages come from.  As a result, I was forced to turn to the Author for clarification. During a long series of meetings, a new text was developed. In doing so, it turned out that there were gaps in many places, and some local errors. The text, moreover, required not only additions, but also abbreviations, as the general editorial awkwardness caused lengthiness, e.g. theorem 2 was proved first in the special case and then only in general, although the general proof does not use the special. I stress that I acted only as a critic. The author was able to clarify each of my doubts, and the revised text which will appear in Colloquium Math., is a completely independent work by Dr. Rokowska. However, it cannot be considered only the result of editorial corrections. It differs from the original text also in terms of content, and the original text, sent to the Faculty Council, was not suitable for printing in my opinion. 

In the printed parts III and II of the dissertation, one can also find ambiguities and omissions, for example, in part III, the second part of the theorem is presented without proof and without any mention of where it comes from. The author provided explanations orally. 

The author has included an interestingly written outline of the history of research from Steiner's theory of systems. It testifies to her reading in this area. In it, the Author mentions her own results very briefly. 

After her doctorate and before her postdoctoral dissertation, Ms. Rokowska published four research papers. Three of them deal with Steiner systems, and one jointly with Z.~Szmerliński relates to technology. The first of these four papers: \emph{Some new constructions of 4-tuples systems} (\CM\ 17, 1967) contains new constructions of $\text{C}(n, 4, 3, 1)$ configurations, i.e., fours in which each three occurs exactly once. Such fours are constructed by the author for certain cardinality n by a different method than Hanani did for those powers earlier, and obtained configurations non-isomorphic with those. In this work, as in later works, one can see the Author's great familiarity with the subject and ability to create complex constructions. In the paper \emph{Some remarks on the number of different triple systems} (Coll. Math.~22, 1971), an estimate of the number of such triple systems for $n\equiv 9 \mod 18$ is given, depending on analogous estimates for certain divisors of the number $n$. 

Attempting a general assessment of Dr. Rokowska's mathematical creativity and achievements is not easy. There is no doubt that her results in Steiner systems theory count and have been noted by specialists, and it would be difficult to surpass her in the structural methods of this theory. However, two serious objections arise.
\begin{enumerate}\itemsep1pt \parskip2pt \parsep1pt
\item Ms. Rokowska's work is extremely limited in terms of subject matter, it is a play on a single string. \emph{Finitist} mathematics, and combinatorics in particular, seems to be experiencing a renaissance. However, Steiner's theory of systems is only a minor part of modern combinatorics. With all caution in assessing the scientific value of this or that issue (any such assessment has something of a divination about it), I fear the sterility of this research. One can easily imagine further estimations of the number of Steiner configurations of increasingly complex type, estimations perhaps of increasingly sharper type, along the path of long and arduous work, but there is no goal in sight at the end of this path, no hypothesis to be proved or disproved. One does not see that one is aiming for an estimate that is in some sense optimal, if only in the case of triples, e.g. to obtain some transparent asymptotics. The methods used before Dr. Rokowska are purely elementary. Perhaps in Steiner's theory of systems (not in all combinatorics certainly) others are not. However, the fact that the candidate does not go beyond completely elementary methods in any of her work, combined with the poor subject matter of her papers, forces one to conclude that her mathematical horizon is tighter than should be required of a postdoctoral fellow. 

\item A postdoctoral fellow is expected to have the ability to properly express mathematical thought and write papers independently. Ms. Dr. Rokowska does not have this skill. This lack is evident when reading Part I of the dissertation, which I have already written about extensively, and the previous manuscripts submitted by the candidate to the Colloquium Mathematicum were also unsuitable for publication and could be published only after alterations, which were made by members of the editorial committee in the course of painstaking work with the Author. Competent reviewers usually cooperating with the editorial board mostly refused to check her papers due to their poor readability.   
\end{enumerate}

I believe that the two allegations I have made here are in some relation to each other and should be treated together. Both of them cause me not to put forward a motion to grant Dr. Barbara Rokowska a postdoctoral degree.

The very long delay in sending this opinion was due to special difficulties in carrying out my task. For this delay I apologize to the Faculty Council.    

\section{\label{BRappBIBD}On the block system designs.} A~block design is an incidence structure consisting of a set $V$ together with a family of subsets $\cB$, called blocks, chosen so that the frequency of the elements satisfies certain conditions that make the collection of blocks exhibit symmetry (balance). We can say that on $V$ the structure is spanned. The distinguished subsets are of identical power, say they have k elements. However, they are not all such subsets and hence incomplete in name. The condition that causes the selected k-element subsets to form a structure is related to the inclusion of selected subsets of lower power in a fixed proportion - in the simplest, most common, and used structure, these are two-element subsets (pairs), and we demand that they occur exactly once. 

Here are the basic definitions related to the study of block systems and their isomorphisms. Without further specification, the term block design usually refers to a balanced incomplete block design (\BIBD\label{BRbibd}), specifically (and also synonymously) a 2-design, which has historically been the most intensively studied type due to its application in experimental design. 

Its generalization is known as a t-design. A design is said to be balanced (up to $t$) if all $t$ subsets of the original set occur in the same number (i.e. $\lambda$) of blocks. Let us collect the list of parameters of a typical \BIBD (=$\BIBD(v,k,t,\lambda$)):
\begin{itemize}\itemsep1pt \parskip2pt \parsep1pt
    \item[v] -- the cardinality of $V$;
    \item[k] -- the size of the blocks; 
    \item[t] -- the size of the subsets to be repeated in the family of blocks;
    \item[$\lambda$] -- the repeating number of $t$ subsets in the design.
\end{itemize}
When $t=k-1$ the block design is called the Steiner block design, usually denoted $\SD(v,k,\lambda))$, and $\SD(v,k)$ when $\lambda=1$ by default.

\begin{definition}[STS]
The Steiner Triple System ($\STS=\STS(n)=\BIBD(n,3,2,1)$) is a $3$-uniform hypergraph on $n$ vertices such that each pair belongs to exactly one edge. 
\end{definition}
\begin{definition}
A design with parameters $t$, $k$, $n$, written $\SD(n,k,t)$, is an $n$-element set $S$ together with a family $\cB$ of $k$-element subsets of $S$ (called blocks) with the property that each $t$-element subset of $S$ is contained in exactly one block. 
\end{definition}
\begin{remark}
In an alternate notation for block designs, an $S(t,k,n)$ would be a $t$-$\BD(n,k,1)$ or  design $\BIBD(n,k,t,1)$. 
\end{remark}

\begin{remark}
In general, one can ask to construct a family $\cB$ of $k$-element subsets of $S$ (called blocks) with the property that each $t$-element subset of $S$ is contained in exactly $\lambda$ block. Such system is denoted $S(t,k,n)$. When $t=k-1$, the system is called Steiner.  
\end{remark}

\begin{definition}[Isomorphism of \STS] Two systems of Stei\-ner tri\-ples, $\cB_1$ and $\cB_2$ on $V$, are isomorphic if there is a permutation of elements (i.e., a bijection transformation $\varphi:V\rightarrow V$) that transforms one system into the other, preserving the relationship between the elements $\varphi(\cB_1)=\cB_2$, where 
\[
\varphi(\cB)=\displaystyle{\bigcup_{\{i_1,i_2,i_3\}\in \cB}}\big\{\{\varphi(i_1),\varphi(i_2),\varphi(i_3)\}\big\}.
\]

If such a transformation exists, we call it an automorphism $(V,\cB_1)$ to $(V,\cB_2)$. Two \STS\ are called non-isomorphic if they are not isomorphic. 
\end{definition}
An \emph{isomorphism} from $\STS(V_1,\cB_1)$ to another  $\STS(V_2,\cB_2)$ is one-to-one map $\varphi$ from $V_1$ to $V_2$ that preserves triples (that is, $\{i_1,i_2,i_3\}\in \cB_1$ if and only if $\{\varphi(i_1),\varphi(i_2),\varphi(i_3)\}\in \cB_2$). An \emph{automorphism} of $\STS(V,\cB)$ is an isomorphism of \STS\ with itself (v.~\cite[p. 2]{Pavo2023:STS13}).

\conflictsofinterest{The author declares that they have no known opposing financial interests or personal relationships that may have influenced the content of this article.
}
 
\funding{This work was funded in part by the Department of Applied Mathematics, Wrocław University of Technology, in the project \textbf{8211204601 MPK: 9130740000}.
}

\acknowledgments{The text was written thanks to the kindness the authors experienced in collecting the information needed for this paper. The primary sources are materials collected at AUWr and APWr. I would like to thank Ms. Anna Rubin-Sieradzka (APWr), as well as colleagues at the Department of Mathematics, Wrocław University of Technology, for their help in reaching the scattered information about Professor Barbara Rokowska's doctoral students.  
}
\begin{center}{\large\bf Resources}\end{center}


\begin{thebibliography}{29}\itemsep3pt \parskip2pt \parsep2pt
\providecommand{\natexlab}[1]{#1}
\providecommand{\url}[1]{\texttt{#1}}
\expandafter\ifx\csname urlstyle\endcsname\relax
  \providecommand{\doi}[1]{doi: #1}\else
  \providecommand{\doi}{doi: \begingroup \urlstyle{rm}\Url}\fi

\bibitem[Alekseev(1974)]{Alekseev1974:STS}
V.~E. Alekseev.
\newblock The number of {S}teiner triple systems.
\newblock \emph{Mat. Zametki}, 15:\penalty0 769--774, 1974.
\newblock ISSN 0025-567X.
\newblock \MR{360304}, \ZBL{0304.05004}, \ZBL{0291.05006}.

\bibitem[{Archiwum \UWr}(1966)]{Rokowska1966}
{Archiwum \UWr}.
\newblock {Teczka przewodu doktorskiego Barbary Rokowskiej}, 9th of March 1966.
\newblock
  \href{https://archiw.uwr.edu.pl/wp-content/uploads/sites/114/2023/01/DR_mat_fiz_chem_1945-2015.pdf?x49020}{Rozprawa
  doktorska 898}. Syg. Akt 34/500.

\bibitem[Arszy{\'{n}}ska(1987)]{Arszynska1987:PhD}
J.~Arszy{\'{n}}ska.
\newblock \emph{{Pewne w{\l}asno{\'{s}}ci symetrycznych system{\'{o}}w i
  kod{\'{o}}w liniowych}}.
\newblock PhD thesis, Instytut Matematyki, Politechnika Wroc{\l}awska,
  Wroc{\l}aw, 1987.

\bibitem[{Baza Danych OPI}(2012)]{BarRok2012:Ludzie}
{Baza Danych OPI}.
\newblock {Barbara Rokowska}, 2012.
\newblock
  \href{https://nauka-polska.pl/\#/profile/scientist?id=60418\&\_k=c9tlpj}{Ludzie
  Nauki ID:60418}, data dostępu=2022-08-23.

\bibitem[Bąk(1991)]{Bak1991:PhD}
M.~Bąk.
\newblock
  \emph{\href{https://omnis-pwr.primo.exlibrisgroup.com/permalink/48OMNIS_TUR/d7ok8p/alma990000687870107668}{O
  rozdzielnych i prawie rozdzielnych systemach bloków}}.
\newblock PhD thesis, Instytut Matematyki i Fizyki Teoretycznej, Politechnika
  Wroc{\l}awska, Wrocław, 1991.

\bibitem[Doyen(1973)]{Doy1973:Steiner}
J.~Doyen.
\newblock Recent results on {S}teiner triple systems.
\newblock In \emph{Finite geometric structures and their applications ({C}entro
  {I}nternaz. {M}at. {E}stivo ({C}.{I}.{M}.{E}.), {II} {C}iclo, {B}ressanone,
  1972)}, Centro Internazionale Matematico Estivo (C.I.M.E.), pages 201--210.
  Ed. Cremonese, Rome, 1973.
\newblock \MR{345837}.

\bibitem[Doyen(1976)]{Doy1976:Steiner}
J.~Doyen.
\newblock Recent developments in the theory of {S}teiner systems.
\newblock In \emph{Colloquio {I}nternazionale sulle {T}eorie {C}ombinatorie
  ({R}oma, 1973), {T}omo {I}}, pages 277--285. Accad. Naz. Lincei, Rome, 1976.
\newblock \MR{437354}.

\bibitem[Doyen and Rosa(1973)]{DoyRos1973:BibSteiner}
J.~Doyen and A.~Rosa.
\newblock A bibliography and survey of {S}teiner systems.
\newblock \emph{Boll. Un. Mat. Ital. (4)}, 7:\penalty0 392--419, 1973.
\newblock \MR{321751}.

\bibitem[Doyen and Rosa(1978)]{DoyRos1978:BibSteiner}
J.~Doyen and A.~Rosa.
\newblock An extended bibliography and survey of {S}teiner systems.
\newblock In \emph{Proceedings of the {S}eventh {M}anitoba {C}onference on
  {N}umerical {M}athematics and {C}omputing ({U}niv. {M}anitoba, {W}innipeg,
  {M}an., 1977)}, Congress. Numer., XX, pages 297--361. Utilitas Math.,
  Winnipeg, MB, 1978.
\newblock \MR{535016}.

\bibitem[Doyen and Rosa(1980)]{DoyRos1980:BibSteiner}
J.~Doyen and A.~Rosa.
\newblock An updated bibliography and survey of {S}teiner systems.
\newblock \emph{Ann. Discrete Math.}, 7:\penalty0 317--349, 1980.
\newblock Topics on Steiner systems. \MR{584420}.

\bibitem[Duda et~al.(2007)Duda, Grząślewicz, Lenczewski, Romanowicz, and
  Weron]{Chmielewski2007:WSA}
R.~Duda, R.~Grząślewicz, R.~Lenczewski, Z.~Romanowicz, and A.~Weron.
\newblock Wrocławska szkoła matematyczna.
\newblock In A.~Chmielewski, editor,
  \emph{\href{http://fbc.pionier.net.pl/id/oai:dbc.wroc.pl:1789}{Wrocławskie
  Środowisko Akademickie : twórcy i ich uczniowie 1945-2005}}, pages
  328--344. Zakład Narodowy im. Ossolińskich, Wrocław, 2007.
\newblock ISBN 978-83-04-04823-2.
\newblock Dodatek do książki:
  \href{http://www.dbc.wroc.pl/publication/1608}{płyta CD}. Zawiera biogramy
  uczonych oraz dokumentację fotograficzną wystawy "Wrocławskie Środowisko
  Akademickie. Twórcy i ich uczniowie", zorganizowanej na Politechnice
  Wrocławskiej w dniach 12-18 listopada 2003 roku i towarzyszącej seminarium
  pod tym samym tytułem.

\bibitem[Erd\H{o}s(1960)]{Erd1960:MR}
P.~Erd\H{o}s.
\newblock {Review \MR{0117188} (22 \#7970)}, 1960.
\newblock ISSN 0025-5629.
\newblock v.~\MR{117188}. Corrected in 2005.

\bibitem[Erd\H{o}s(1967)]{Hanani1967:MR}
P.~Erd\H{o}s.
\newblock {Review \MR{0213253 (35 \#4117)}}, 1967.
\newblock ISSN 0025-5629.
\newblock v.~\MR{0213253}.

\bibitem[Erd{\H{o}}s(1960)]{Erdos1960:Zahorski}
P.~Erd{\H{o}}s.
\newblock About an estimation problem of {Zahorski}.
\newblock \emph{Colloq. Math.}, 7:\penalty0 167--170, 1960.
\newblock ISSN 0010-1354.
\newblock \doi{10.4064/cm-7-2-167-170}.
\newblock \ZBL{0106.27701}.

\bibitem[Hanani(1960)]{Hanani1960:Quadruple}
H.~Hanani.
\newblock On quadruple systems.
\newblock \emph{Canadian J. Math.}, 12:\penalty0 145--157, 1960.
\newblock ISSN 0008-414X.
\newblock \doi{10.4153/CJM-1960-013-3}.
\newblock \MR{111696}.

\bibitem[Hanani(1963)]{Hanani1963:Tact}
H.~Hanani.
\newblock On some tactical configurations.
\newblock \emph{Canadian J. Math.}, 15:\penalty0 702--722, 1963.
\newblock ISSN 0008-414X.
\newblock \doi{10.4153/CJM-1963-069-5}.
\newblock \MR{157908}.

\bibitem[Kwapień(2016)]{Kwapien2015:CRN}
S.~Kwapień.
\newblock Czesław ryll-nardzewski (7 x 1926 – 18 ix 2015).
\newblock \emph{Rocznik Polskiej Akademii Umiejętności}, 2015/2016:\penalty0
  \href{https://pau.krakow.pl/Rocznik_PAU/2015_2016/Rocznik_2015_2016_zmarli_Ryll_Nardzewski.pdf}{176--178},
  2016.

\bibitem[Pavone(2023)]{Pavo2023:STS13}
M.~Pavone.
\newblock A visual representation of the {Steiner} triple systems of order 13.
\newblock \emph{Art Discrete Appl. Math.}, 6\penalty0 (3):\penalty0 22, 2023.
\newblock ISSN 2590-9770.
\newblock \doi{10.26493/2590-9770.1564.2b8}.
\newblock Id/No p3.04. \ZBL{1509.05036}.

\bibitem[Pawlik(1984)]{Paw1984:PhD}
B.~Pawlik.
\newblock \emph{{Prawie rozdzielne systemy blokow}}.
\newblock Rozprawa doktorska, Instytut Matematyki, Politechnika Wroc{\l}awska,
  Wroc{\l}aw, 10.07 1984.

\bibitem[Piegatowski(1993)]{Piegatowski1993:PhD}
D.~Piegatowski.
\newblock \emph{{Łuki zupełne w systemach trójek Steinera}}.
\newblock PhD thesis, Instytut Matematyki i Fizyki Teoretycznej, Politechnika
  Wroc{\l}awska, Wroc{\l}aw, 1993.

\bibitem[Redaktor\_WM(2017)]{WM2316:PTM}
N.~Redaktor\_WM.
\newblock Wykaz członków ptm.
\newblock \emph{Wiadomości Matematyczne}, 8\penalty0 (3):\penalty0 189--201,
  2017.
\newblock ISSN 2543-991X.
\newblock \doi{10.14708/wm.v8i3.2316}.

\bibitem[Rokowska(1966)]{BRok1966:PhD}
B.~Rokowska.
\newblock \emph{Pewne nowe konstrukcje systemów Steinera}.
\newblock Rozprawa doktorska, Wroclaw University. Faculty of Mathematics,
  Physics and Chemistry, Wrocław, 9th of March 1966.
\newblock \emph{Certain new constructions Steiner's systems}. Archiwum
  Uniwersytetu Wrocławskiego
  \href{https://archiw.uwr.edu.pl/wp-content/uploads/sites/114/2023/01/DR_mat_fiz_chem_1945-2015.pdf?x49020}{Rozprawa
  doktorska 898} Syg. Akt 34/500.

\bibitem[Rokowska(1976)]{Rokowska1976}
B.~Rokowska.
\newblock \emph{Nieizomorficzne systemy trójek Steinera}.
\newblock Rozprawa habilitacyjna, Wroclaw University. Faculty of Mathematics,
  Physics and Chemistry, Wrocław, 15th of December 1976.
\newblock \emph{Non-isomorphic Steiner's triples systems}. Archiwum
  Uniwersytetu Wrocławskiego
  \href{https://archiw.uwr.edu.pl/wp-content/uploads/sites/114/2023/01/hab_matfizchem_1945-2015.pdf?x49020}{Rozprawa
  habilitacyjna 283} Syg. Akt 35/147.

\bibitem[Rokowska(1996)]{Rok1996:Jezus}
B.~Rokowska.
\newblock \emph{Jezus we Wrocławiu}.
\newblock Poligrafia Inspektoratu Towarzystwa Salezjańskiego, ul. Konfederacka
  6, 30-306 Kraków, 1996.
\newblock ISBN 83-86473-14-2.

\bibitem[Rokowska(2005)]{Rok2005:Jezus}
B.~Rokowska.
\newblock \emph{Wyzwolenie ku radości}.
\newblock Wydawnictwo św. Antoniego, ul. Konfederacka 6, 30-306 Kraków, 2005.
\newblock ISBN 83-88598-62-7.

\bibitem[Tyszkowski(2000)]{Tysz2000:Ossoliniana}
W.~Tyszkowski.
\newblock Zarys historii działu starych druków zakładu narodowego im.
  ossolińskich 1918--1998.
\newblock \emph{Ossoliniana. Czasopismo ZNiO},
  \href{https://ossolineum.pl/index.php/czasopismo-znio/numery-archiwalne/zeszyt-11-2000/}{11}:\penalty0
  217--226, 2000.
\newblock ISSN 1230-221X.
\newblock Czasopismo Zakładu Narodowego im. Ossolińskich. Dział:
  Ossoliniana.

\bibitem[{University of Wroclaw Archives}(1959)]{BRokowska1954:UWr}
{University of Wroclaw Archives}.
\newblock Personal file of a student.
\newblock University of Wroclaw Archives; File reference number 154, Barbara
  Rokowska, 1959.
\newblock Period 1.08.1954--1859.

\bibitem[Wilczyńska(1981)]{Wil1981:PhD}
K.~Wilczyńska.
\newblock \emph{Rozdzielne systemy czwórek $\text{RB}(v,4,3)$}.
\newblock Rozprawa doktorska, Instytut Matematyki i Fizyki Teoretycznej,
  Politechnika Wrocławska, Wrocław, 24.06 1981.
\newblock \emph{A construction of resolvable quadruple systems}.
  \href{https://dona.pwr.edu.pl/szukaj/default.aspx?nrewid=058968}{Raporty
  Inst. Mat. PWr. 1981, Ser. PRE, nr 2, 38p. I18/1981/P-002}.

\bibitem[{Wroclaw University of Technology Archives}(1959)]{BRokowska1959:PWr}
{Wroclaw University of Technology Archives}.
\newblock Personal file of an employee.
\newblock Wroclaw University of Technology Archives; File reference number
  1431/176(36-E-205), Barbara Rokowska, 1959.
\newblock Period 1.08.1959--2002.

\end{thebibliography}
\setcounter{section}{0}
\selectlanguage{polish}
\Polskitrue
\subjclass[2010]{05B05; 01A50; 01A55; 01A60}

\keywords{Systemy Steinera; system bloków; liczby pierwsze; ciała Galois; konstrukcje kombinatoryczne}

\begin{center}
{\bf 
Badania kombinatoryczne Barbary Rokowskiej\\ z jej obszerną biografią (1926--2012)\\
Krzysztof J. Szajowski}
\end{center}
\vspace{-2ex}
\begin{abstract}
Omawiamy znaczenie kilku interesujących wyników Barbary Rokowskiej o konstrukcjach kombinatorycznych. Jej zainteresowanie matematyką skończoną i teorią liczb rozpoczęło się od uścislenia i uszczegółowienia pewnej pracy Erdősa. Następnie Rokowska i Schinzel rozwiązali problem postawiony przez Paula Erdősa dotyczący istnienia liczb pierwszych pewnego rodzaju. Jej kolejne prace uwypukliły trudność w konstrukcjach systemów Steinera o pewnych własnościach i pokazała znaczenie rygorystycznych technik dowodowych w tym obszarze matematyki. Jest to pierwsze takie zestawienie głównych wyników uzyskanych przez Rokowską, jej współpracowników i doktorantów.

Biografia Barbary Rokowskiej została dodana jako dodatek.
\end{abstract}
\setcounter{section}{0}
\selectlanguage{english}
\Polskifalse

\vspace{-5ex}

\Koniec
\end{document}